\documentclass[sn-mathphys,Numbered]{sn-jnl}

\usepackage{graphicx}%
\usepackage{multirow}%
\usepackage{amsmath,amssymb,amsfonts}%
\usepackage{amsthm}%
\usepackage{mathrsfs}%
\usepackage[title]{appendix}%
\usepackage{xcolor}%
\usepackage{textcomp}%
\usepackage{manyfoot}%
\usepackage{booktabs}%
\usepackage{algorithm}%
\usepackage{algorithmicx}%
\usepackage{algpseudocode}%
\usepackage{listings}%
\usepackage{tabularx}
\usepackage[all]{xy}
\usepackage{tikz}
\usetikzlibrary{shapes,decorations}
\usetikzlibrary {arrows.meta}
\usepackage[textsize=scriptsize,backgroundcolor=white]{todonotes}

\theoremstyle{plain}%
\newtheorem{theorem}{Theorem}
\newtheorem{proposition}[theorem]{Proposition}%
\newtheorem{lemma}[theorem]{Lemma}
\newtheorem{corollary}[theorem]{Corollary}

\theoremstyle{remark}%
\newtheorem{example}{Example}%
\newtheorem{remark}{Remark}%

\theoremstyle{definition}%
\newtheorem{definition}{Definition}%

\def\mcl#1{\mathcal{#1}}
\def\bracket#1{\left\langle #1\right\rangle}
\def\hil{\mcl{H}}
\def\nn{\nonumber}
\def\opn{\operatorname}
\def\mr{\mathrm}

\def\alg{\mcl{A}}
\def\modu{\mcl{M}}
\def\moduu{\mcl{N}}

\def\Bbracket#1{\bigg\langle #1\bigg\rangle}
\def\sbracket#1{\langle #1\rangle}

\def\zhat{\hat{\mcl{Z}}}

\def\ly{L^2(\mcl{Y})}
\def\lz{L^2(\mcl{Z})}
\def\lx{L^2(\mcl{X})}

\DeclareSymbolFont{EulerExtension}{U}{euex}{m}{n}
\DeclareMathSymbol{\euintop}{\mathop} {EulerExtension}{"52}
\DeclareMathSymbol{\euointop}{\mathop} {EulerExtension}{"48}

\DeclareMathOperator*{\esssup}{ess\,sup}

\begin{document}

\title[]{Koopman spectral analysis of skew-product dynamics on Hilbert $C^*$-modules}


\author[1]{\fnm{Dimitrios} \sur{Giannakis}}\email{dimitrios.giannakis@dartmouth.edu}

\author*[2,3]{\fnm{Yuka} \sur{Hashimoto}}\email{yuka.hashimoto@ntt.com}

\author[3,4]{\fnm{Masahiro} \sur{Ikeda}}\email{masahiro.ikeda@riken.jp}

\author[5,3]{\fnm{Isao} \sur{Ishikawa}}\email{ishikawa.isao.zx@ehime-u.ac.jp}

\author[1]{\fnm{Joanna} \sur{Slawinska}}\email{joanna.m.slawinska@dartmouth.edu}

\affil[1]{\orgdiv{Department of Mathematics}, \orgname{Dartmouth University}, \orgaddress{\street{27 N. Main Street}, \city{Hanover}, \postcode{03755}, \state{New Hampshire}, \country{USA}}}

\affil*[2]{\orgdiv{NTT Network Service Systems Laboratories}, \orgname{NTT Corporation}, \orgaddress{\street{3-9-11, Midori-cho}, \city{Musashino}, \postcode{180-8585}, \state{Tokyo}, \country{Japan}}}

\affil[3]{\orgdiv{Center for Advanced Intelligence Project}, \orgname{RIKEN}, \orgaddress{\street{1-4-1 Nihonbashi}, \city{Chuo}, \postcode{103-0027}, \state{Tokyo}, \country{Japan}}}

\affil[4]{\orgdiv{Faculty of Science and Technology}, \orgname{Keio University}, \orgaddress{\street{3-14-1, Hiyoshi}, \city{Kohoku, Yokohama}, \postcode{223-8522}, \state{Kanagawa}, \country{Japan}}}

\affil[5]{\orgdiv{Center for Data Science}, \orgname{Ehime University} \orgaddress{\street{2-5, Bunkyo-cho}, \city{Matsuyama}, \postcode{790-8577}, \state{Ehime}, \country{Japan}}}


\abstract{We introduce a linear operator on a Hilbert $C^*$-module for analyzing skew-product dynamical systems.
The operator is defined by composition and multiplication.
We show that it admits a decomposition in the Hilbert $C^*$-module, called eigenoperator decomposition, that generalizes the concept of the eigenvalue decomposition.
This decomposition reconstructs the Koopman operator of the system in a manner that represents the continuous spectrum through eigenoperators. 
In addition, it is related to the notions of cocycle and Oseledets subspaces and it is useful for characterizing coherent structures under skew-product dynamics. We present numerical applications to simple systems on two-dimensional domains.}




\maketitle
\section{Introduction}

\subsection{Background and motivation}
Operator-theoretic methods have been used extensively in analysis and computational techniques for dynamical systems.  
Let $f:\mcl{X}\to\mcl{X}$ be a dynamical system on a state space $\mathcal X$.
Then, the Koopman operator $U_f$ associated with $f$ is defined as a composition operator on an $f$-invariant function space $\mcl{F}$ on $\mcl{X}$,
\begin{equation*}
    U_fv=v\circ f, 
\end{equation*}
for $ v \in \mathcal F$ \cite{koopman31,KoopmanVonNeumann32}. In many cases, $\mcl{F}$ is chosen as a Banach space or Hilbert space, such as the Lebesgue spaces $L^p(\mcl{X})$ for a measure space $\mcl{X}$ and the Hardy space $H^p(\mathbb{D})$ on the unit disk $\mathbb{D}$, where $p\in [1,\infty]$. Meanwhile, the Perron--Frobenius, or transfer, operator associated with $f$ is defined as the adjoint $P_f$ of the Koopman operator acting on the continuous dual $\mathcal F'$ of $\mathcal F$, i.e., $P_f \nu = \nu \circ U_f$ for $\nu \in \mathcal F'$ \cite{Baladi00}. In a number of important cases (e.g., $\mathcal F = L^p(\mathcal X)$ with $ p \in [1, \infty)$ or $\mathcal F = C(\mathcal X)$ for a compact Hausdorff space $\mathcal X$), $\mathcal F'$ can be identified with a space of measures on $\mathcal X$; the transfer operator is then identified with the pushforward map on measures, $P_f \nu = \nu \circ f^{-1}$. When $\mathcal F$ has a predual, $\mathcal F_* \subseteq \mathcal F'$ it is common to define $P_f$ as the predual of the Koopman operator, i.e., $ (U_f v) \nu = v(P_f \nu)$; an important such example is $\mathcal F = L^\infty(\mathcal X)$ with $\mathcal F_* = L^1(\mathcal X)$. A central tenet of modern ergodic theory is to leverage the duality relationships between $f: \mathcal X \to \mathcal X$, $U_f : \mathcal F \to \mathcal F$, and $P_f : \mathcal F' \to \mathcal F'$ to characterize properties of nonlinear dynamics such ergodicity, mixing, and existence of factor maps, using linear operator-theoretic techniques \cite{EisnerEtAl15}. 

Starting from work in the late 1990s and early 2000s \cite{Froyland97,DellnitzJunge99,DellnitzEtAl00,Mezic05}, operator-theoretic techniques have also proven highly successful in data-driven applications \cite{GiannakisEtAl15,klus17,ishikawa18,hashimoto20}.
A primary such application is the modal decomposition (e.g., \cite{Schmid10,RowleyEtAl09,williams15,kawahara16,hassan17,RosenfeldEtAl22}).
This approach applies eigenvalue decomposition to the Koopman operator to identify the long-term behavior of the dynamical system.
Assume $\mcl{F}$ is a Hilbert space equipped with an inner product $\bracket{\cdot,\cdot}$, and $U_f$ is normal, bounded, and diagonalizable, with eigenvalues $\lambda_1,\lambda_2,\ldots \in \mathbb C$ and corresponding basis of orthonormal eigenvectors $v_1,v_2,\ldots \in \mathcal F$.
Then, for $u\in\mcl{F}$ and a.e.\ $x \in \mathcal X$ we have $u(f^i(x)) =U_f^iu(x)=\sum_{j=1}^{\infty}\lambda_j^iv_j(x)\bracket{v_j,u}$, where $i \in \mathbb N$ represents discrete time. 
Therefore, the time evolution of observables is described by the Koopman eigenvalues and corresponding eigenvectors.
By computing the eigenvalues of the Koopman operator, we obtain oscillating elements and decaying elements in the dynamical system.

Several attempts have been made to generalize the above decomposition to the case where the Koopman operator has continuous or residual spectrum.
Korda et al.~\cite{korda20} approximate the spectral measure of the Koopman operator on $L^2(\mathcal X)$ for measure-preserving dynamics using Christoffel--Darboux kernels in spectral space.
Slipantschuk et al.~\cite{slipantschuk20} consider a riddged Hilbert space and extend the Koopman operator to a space of distributions so that it becomes compact. Colbrook and Townsend \cite{ColbrookTownsend21} employ a residual-based approach that consistently approximates the spectral measure by removing spurious eigenvalues from DMD-type spectral computations.
Spectrally approximating the Koopman operator in measure-preserving, ergodic flows by compact operators on reproducing kernel Hilbert spaces (RKHSs) has also been investigated \cite{das21}. However, dealing with continuous and residual Koopman spectra is still a challenging problem. 

On the transfer operator side, popular approximation techniques are based on the Ulam method \cite{Ulam64}. The Ulam method has been shown to yield spectrally consistent approximations for particular classes of systems such as expanding maps and Anosov diffeomorphisms on compact manifolds \cite{Froyland97}. In some cases, spectral computations from the Ulam method has been shown to recover eigenvalues of transfer operators on anisotropic Banach spaces adapted to the expanding/contracting subspaces of such systems \cite{BlankEtAl02}; however, these results depend on carefully chosen state space partitions that may be hard to construct in high dimensions and/or under unknown dynamics. Various modifications of the basic Ulam method have been proposed that are appropriate for high-dimensional applications; e.g., sparse grid techniques \cite{JungeKoltai09}.

\subsection{Skew-product dynamical systems}

We focus on measure-preserving skew-product systems in discrete time, $T(y,z)=(h(y),g(y,z))$, or continuous-time, $\Phi_t(y,z)=(h_t(y),g_t(y,z))$, on a product space $ \mathcal X = \mcl{Y}\times\mcl{Z}$. Here, $\mcl{Y}$ and $\mcl{Z}$ are measure spaces, oftentimes referred to as the ``base'' and ``fiber'', respectively.
In such systems, the driving dynamics on $\mcl{Y}$ is autonomous, but the  dynamics on $\mcl{Z}$ depends on the configuration $ y \in \mcl{Y}$. In many cases, one is interested in the time-dependent fiber dynamics, rather than the autonomous dynamics on the base.
A typical example of skew-product dynamics is Lagrangian tracer advection under a time-dependent fluid flow \cite{froyland10b,froyland17,giannakis20}, where $\mathcal Y$ is the state space of the fluid dynamical equations of motion and $\mathcal Z$ is the spatial domain where tracer advection takes place.

A well-studied approach for analysis of skew-product systems involves replacing the spectral decomposition of Koopman/transfer operators acting on functions on $\mathcal X$ by decomposition of associated operator \emph{cocycles} acting on functions on $\mathcal Y$ using multiplicative ergodic theorems. In a standard formulation of the multiplicative ergodic theorem, first proved by Oseledets \cite{Oseledets68}, one considers an invertible measure-preserving map $h:\mcl{Y}\to\mcl{Y}$ and the cocycle generated by a matrix-valued map $A$ on $\mcl{Y}$. 
The multiplicative ergodic theorem then shows the existence of subspaces $\mcl{V}_1(y),\ldots,\mcl{V}_k(y)$ such that $A(y)\mcl{V}_j(y)=\mcl{V}_j(h(y))$.
The subspace $\mcl{V}_j(y)$ is called an Oseledets subspace (or equivariant subspace). Each Oseledets subspace has an associated Lyapunov exponent and associated covariant vectors, which are the analogs of the eigenvalues and eigenvectors of Koopman/transfer operators, respectively, in the setting of cocycles. Since its inception, the multiplicative ergodic theorem has been extended in many ways to infinite-dimensional operator cocycles \cite{Ruelle68,thieullen87,Schaumloffel91,FroylandEtAl10c,gonzalez14}. Under appropriate quasi-compactness assumptions, it has been shown that the Lyapunov exponent spectrum is at most countably infinite and the associated Oseledets subspaces are finite-dimensional, e.g., \cite[Theorem~A]{gonzalez14}. 

A primary application of Oseledets decompositions is the detection of coherent sets and coherent structures in natural and engineered systems~\cite{FROYLAND10a}.
A family of sets $\{\mcl{S}(y)\}_{y\in\mcl{Y}}$ is called coherent if $\nu(\mcl{S}(y)\bigcap g(h(y),\mcl{S}(y)))/\nu(\mcl{S}(y))$ is large for a reference measure $\nu$ on $\mcl{Z}$.
If an Oseledets subspace $\mcl{V}(y)$ with respect to the transfer operator cocycle $U_{g(y,\cdot)}^{-1}$ is represented as $\mcl{V}(y)=\opn{Span}\{v_y\}$ for a covariant vector $v_y\in\lz$ satisfying $U_{g(y,\cdot)}v_{h(y)}=v_y$,
then setting $\mcl{S}(y)$ to a level set of $v_y$ leads to a family of coherent sets.
Finite-time coherent sets and Lagrangian coherent structures as the boundaries of the finite-time coherent sets have also been studied~\cite{froyland10b,froyland18,froyland15}.

\subsection{Eigenoperator decomposition}

In this paper, we investigate a different approach to deal with continuous and residual spectra of Koopman operators on $L^2$ associated with skew-product dynamical systems. We propose a new decomposition, called {\em eigenoperator decomposition}, which reconstructs the Koopman operator from multiplication operators acting on certain subspaces, referred to here as generalized Oseledets spaces. These multiplication operators are obtained by solving an eigenvalue-type equation, but they can individually have continuous spectrum. Intuitively, this decomposition provides a factorization of the (potentially continuous) spectrum of the underlying Koopman operator into the spectra of eigenoperator families. 

Our approach is based on the theory of Hilbert $C^*$-modules \cite{lance95}, which generalizes Hilbert space theory by replacing the complex-valued inner product by a product that takes values in a $C^*$-algebra.
In this work, we employ the $C^*$-algebra of bounded linear operators on $L^2(\mcl{Z})$, denoted by $\mcl{B}(\lz)$.
A standard operator-theoretic approach for skew-product dynamics is to define the Koopman or transfer operator on the product Hilbert space $\mathcal H = L^2(\mcl{Y})\otimes L^2(\mcl{Z})$ \cite{giannakis20,FroylandKoltai23}. 
In contrast, here we consider the Hilbert $C^*$-module $\mathcal M = \ly\otimes\mcl{B}(\lz)$ over $\mcl{B}(\lz)$.
By considering $\mcl{B}(\lz)$ instead of $\lz$, we aim to push information about the continuous spectrum of the Koopman operator onto the $C^*$-algebra $\mcl{B}(\lz)$. 

\begin{figure}[t]
\begin{tikzpicture}
\node at (-5,5)[below,align=left] {$\Phi_t(y,z)=(h_t(y),g_t(y,z))$\\ \quad: Skew product flow};
\draw[-latex, ultra thick] (-2.8,4.5) to (-2.1,4.5);
\node at (0.5,5)[below,align=left] {$K_{\Phi_t}$: Linear operator on \\ \phantom{$K_{\Phi_t}$: }the Hilbert $C^*$-module\\ \phantom{$K_{\Phi_t}$: }$\modu=\ly\otimes \mcl{B}(\lz)$};
\draw[-latex, ultra thick] (3,4.5) to (3.7,4.5);
\node at (5.2,5)[below,align=left] {$L_{\Phi}$: Generator\\ \phantom{$L_{\Phi}$: }of $K_{\Phi}$};
\draw[-latex,ultra thick,blue] (-5,4) to (-5,3.3);
\node at (-5,3.3)[below,align=left] {$U_{\Phi_t}$: Koopman operator on \\ \phantom{$U_{\Phi_t}$: }a Banach space \\ \phantom{$U_{\Phi_t}$: }$\moduu=C(\mcl{Y})\otimes \lz$};
\draw[-latex,ultra thick,blue] (-5,1.8) to (-5,0.3);
\node at (-5,0.3)[below,align=left] {$V_{\Phi}$: Generator\\ \phantom{$V_{\Phi}$: }of $U_{\Phi_t}$};
\draw[-latex,ultra thick] (5.2,4) to (5.2,3.1);
\node at (-1,3.3)[below,blue] {Eigenvector};
\node at (2,3.3)[below,red] {Eigenoperator};
\draw[-latex,ultra thick,red] (2,2.9) to (1,2);
\draw[-latex,ultra thick,blue] (-1,2.9) to (-1,2);
\node at (2.3,2.7)[draw,below,align=left] {\underline{Decomposition} (Theorems~\ref{prop:detailed_decom} and \ref{prop:hat_w_conti})\\ $L_{\Phi}\textcolor{blue}{\hat{w}_{s,j}}=\textcolor{blue}{\hat{w}_{s,j}}\cdot\textcolor{red}{\hat{N}_{s,j}}$\quad $\textcolor{blue}{\hat{w}_{s,j}}\in\modu$, $\textcolor{blue}{\hat{w}_{s,j}}=w_sp_j$\\ $(j=1,2,\ldots,\ s\in\mathbb{R})$\ \, $w\cdot\textcolor{red}{\hat{N}_{s,j}}(y)=w(y)(\textcolor{red}{V_{\Phi}}p_j)(h_s(y))$};
\draw[-latex,ultra thick,blue] (-3.7,-0.2) to (-3,-0.2);
\node at (-1.2,0.3)[below,align=left] {$\mcl{V}_{j}$: Invariant\\ \phantom{$\mcl{V}_{j}$: }subspace of $V_{\Phi}$};
\draw[-latex,ultra thick,blue] (0.5,-0.2) to node [above] {Fix $y$} (1.2,-0.2);
\node at (4,0.3)[below,align=left] {$\mcl{V}_{j}(y)$: Oseledets space on $\lz$\\ $p_j(y)$: Projection onto $\mcl{V}_j(y)$};
\draw[-latex,ultra thick,blue] (4,0.3) to (4,1);
\end{tikzpicture}
\caption{Overview of eigenoperator decomposition for continuous-time systems.}
    \label{fig:overview}
\end{figure}
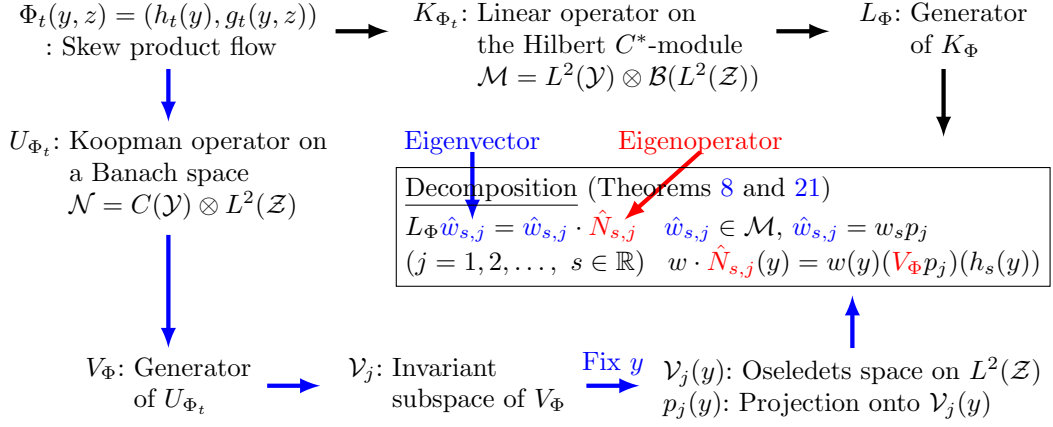

In more detail, starting from discrete-time systems, we define a $\mathcal B(L^2(\mathcal Z))$-linear operator $K_T$ on $\mathcal M$, which can be thought of as a lift of the standard Koopman operator on $\mathcal H$ to the Hilbert $C^*$-module setting. In addition, $K_T$ can be used to reconstruct the Koopman operator of the full skew-product system on $\mathcal X$. We show that $K_T$ admits a decomposition
\begin{equation}
    K_T \hat w_{i,j} = \hat w_{i,j} \cdot \hat M_{i,j},
    \label{eq:decomp}
\end{equation}
where $\hat M_{i,j}$ is a $\mathcal B(L^2(\mathcal Z))$-linear multiplication operator (which we call eigenoperator), and $\hat w_{i,j} \in \mathcal M$ are eigenvectors associated with the operator cocycle on $L^2(\mathcal Z)$ induced by the skew-product dynamics.  We also derive an analogous version of~\eqref{eq:decomp} for continuous-time systems, formulated in terms of the generator of the Koopman group $ \{ K_{\Phi_t} \}_{t \in \mathbb R}$ acting on $\mathcal M$. A schematic overview of our approach for the continuous-time case is displayed in Fig.~\ref{fig:overview}. 

The eigenoperator decomposition~\eqref{eq:decomp} and its continuous-time variant have associated equivariant subspaces of $L^2(\mathcal Z)$ as in the multiplicative ergodic theorem. In particular, to each eigenoperator $\hat M_{i,j}$ there is an associated family $\{\mathcal V_j(y)\}_{y \in \mathcal Y}$ of closed subspaces $\mathcal V_j(y) \subseteq L^2(\mathcal Z)$ such that $U_{g(y,\cdot)}$ maps vectors in $\mathcal V_j(h(y))$ to vectors in $\mathcal V_j(y)$. Since we consider cocycles generated by unitary Koopman/transfer operators, the equivariant subspaces $\mathcal V_j(y)$ can be infinite-dimensional. Therefore, we call them generalized Oseledets subspaces. Spectral analysis of $\hat M_{i,j}$ then reveals coherent structures under the skew-product dynamics. 

The rest of this paper is organized as follows. In Section~\ref{sec:discrete_systems}, we derive our eigenoperator decomposition for discrete-time systems, and establish the correspondence between $K_T$ and the Koopman operator. We illustrate the decomposition in Section~\ref{sec:discrete_time_examples} by means of analytical examples with fiber dynamics on abelian and non-abelian groups. In these examples, the generalized Oseledets subspaces can be constructed explicitly, which provides intuition about the behavior of eigenoperator decomposition. In Section~\ref{sec:continuous_time_systems}, we describe the construction of the infinitesimal generator and the associated eigenoperator decomposition for continuous-time systems. Section~\ref{sec:continuous_time_examples} contains numerical applications of the decomposition for continuous-time systems to simple time-dependent flows in two-dimensional domains. Section~\ref{sec:conclusions} contains a conclusory discussion. The paper includes an Appendix collecting auxiliary results.    


\section{Discrete-time systems}\label{sec:discrete_systems}
\subsection{Skew product system and Koopman operator on Hilbert space}
Let $\mcl{Y}$ and $\mcl{Z}$ be separable measure spaces equipped with measures $\mu$ and $\nu$, respectively
and let $\mcl{X}=\mcl{Y}\times\mcl{Z}$, the direct product measure space of $\mcl{Y}$ and $\mcl{Z}$.
Let $h:\mcl{Y}\to\mcl{Y}$ be a measure preserving and invertible map and let $g:\mcl{X}\to\mcl{Z}$ be a measurable map such that $g(y,\cdot)$ is measure preserving and invertible for any $y\in\mcl{Y}$.
Consider the following skew product transformation $T$ on $\mcl{X}$:
\begin{equation*}
 T(y,z)=(h(y),g(y,z)).
\end{equation*}
We consider the Koopman operator $U_T$ on $L^2(\mcl{X})$.
Note that since $\ly$ and $\lz$ are separable, their tensor product $\ly\otimes\lz$ satisfies
\begin{align*}
\ly\otimes \lz\simeq L^2(\mcl{X}).
\end{align*}
\begin{definition}\label{def:koopman_discrete}
The Koopman operator $U_T$ on $L^2(\mcl{X})$ is defined as
\begin{equation*}
 U_Tf=f\circ T
\end{equation*}
for $f\in L^2(\mcl{X})$.
\end{definition}

Since $T$ is measure preserving, the Koopman operator $U_T$ is an unitary operator, but $U_T$ does not always have an eigenvalue decomposition since it has continuous spectrum in general.

\subsection{Operator on Hilbert $C^*$-module related to the Koopman operator}
We extend the Koopman operator $U_T$ to an operator on a Hilbert $C^*$-module.
We first introduce Hilbert $C^*$-module~\cite{lance95,hashimoto23}.
\begin{definition}
For a module $\modu$ over a $C^*$-algebra $\alg$,
a map $\bracket{\cdot,\cdot}_{\modu}:\modu\times\modu\to\alg$ is referred to as an $\alg$-valued {inner product} if it is $\mathbb{C}$-linear with respect to the second variable and has the following properties:
For $w_1,w_2,w_3\in\modu$ and $a,b\in\alg$, 
\begin{enumerate}
 \item $\bracket{w_1,w_2 a+w_3 b}_{\modu}=\bracket{w_1,w_2}_{\modu}a+\bracket{w_1,w_3}_{\modu}b$,
 \item $\bracket{w_1,w_2}_{\modu}=\bracket{w_2,w_1}_{\modu}^*$,
 \item $\bracket{w_1,w_1}_{\modu}$ is positive,
 \item If $\bracket{w_1,w_1}_{\modu}=0$ then $w_1=0$.
\end{enumerate}
If $\bracket{\cdot,\cdot}$ satisfies the conditions 1$\sim$3, but not 4, then it is called a semi-inner product.
Let $\Vert w\Vert_{\modu} =\Vert\bracket{w,w}_{\modu}\Vert_\alg^{1/2}$ for $w\in\modu$. 
Then $\Vert\cdot\Vert_{\modu}$ is a norm in $\modu$.
\end{definition}
\begin{definition}
A Hilbert $C^*$-module over $\alg$ or Hilbert $\alg$-module is a module over $\alg$ equipped with an $\alg$-valued inner product and complete with respect to the norm induced by the $\alg$-valued inner product.
\end{definition}

Let $\alg$ be the $C^*$-algebra $\mcl{B}(\lz)$.
Let 
\begin{align*}
\mcl{M}=\ly\otimes\alg,
\end{align*}
i.e., the (right) Hilbert $\alg$-module defined by the tensor product of the Hilbert $\mathbb{C}$-module $\ly$ and (right) Hilbert $\alg$-module $\alg$~\cite{lance95}.
We now define an operator on a Hilbert $C^*$-module.
\begin{definition}
We define the a right $\alg$-linear operator $K_T$ on $\modu$ (i.e., $K_T$ is linear and satisfies $K_T(wa)=(K_tw)a$ for all $a\in\alg$ and $w\in\modu$) by
\begin{equation*}
K_T(v\otimes a)(y)=v(h(y))U_{g(y,\cdot)}a
\end{equation*}
for $v\in\ly$, $a\in\alg$, and $y\in\mcl{Y}$.
Here, for $y\in\mcl{Y}$, $U_{g(y,\cdot)}$ is the Koopman operator on $\lz$ with respect to the map $g(y,\cdot)$.
\end{definition}

The well-definedness of $K_T$ is not trivial.
The following proposition shows the well-definedness of $K_T$ as an operator from $\modu$ to $\modu$.
\begin{proposition}\label{prop:K_TinM}
The operator $K_T$ is a right $\alg$-linear unitary operator from $\modu$ to $\modu$.
\end{proposition}
\noindent The proof of Proposition~\ref{prop:K_TinM} is documented in Appendix.

The next proposition shows the relationship between $U_T$ and $K_T$, which enables us to connect existing studies of Koopman operators with our framework.
\begin{proposition}\label{prop:UandK}
Let $\{\gamma_i\}_{i=1}^{\infty}$ be an orthonormal basis of $\lz$.
Let $\iota_i:\lx\to\modu$ and $P_i:\modu\to \lx$ be linear operators defined as $v\otimes u\mapsto v\otimes u\gamma_i'$ and $v\otimes a\mapsto v\otimes a\gamma_i$, respectively.
Then, we have $U_T=P_iK_T\iota_i$ for any $i=1,2,\ldots$.
Moreover, we have $K_T=\sum_{i=1}^{\infty}\iota_i U_T P_i$, where the sum converges strongly to $K_T$ in $\modu$.
\end{proposition}
\begin{proof}
For $v\in\ly$, $u\in\lz$, $y\in\mcl{Y}$, and $i=1,2,\ldots$, we have
\begin{align}
K_T\iota_i(v\otimes u)(y)
&=K_T(v\otimes u\gamma_i')(y)
=v(h(y))U_{g(y,\cdot)}u\gamma_i'\nn\\
&=v(h(y))u(g(y,\cdot))\gamma_i'=\iota_iU_T(v\otimes u).\label{eq:UandK}
\end{align}
By acting $P_i$ on the both sides of Eq.~\eqref{eq:UandK}, we have $P_iK_T\iota_i(v\otimes u)=U_T(v\otimes u)$.
Since $U_T$ and $K_T$ are bounded, $P_iK_T\iota_i=U_T$ holds on $\modu$.
Moreover, for $v\in\ly$ and $a\in\alg$, we have
\begin{align}
\sum_{i=1}^{\infty}\iota_iP_i(v\otimes a)=\sum_{i=1}^{\infty}v\otimes a\gamma_i\gamma_i'=v\otimes a,\label{eq:identity}
\end{align}
where the convergence is the strong convergence.
Since $K_T$ is bounded, by Eqs.~\eqref{eq:UandK} and \eqref{eq:identity}, $K_T(v\otimes a)=\sum_{i=1}^{\infty}\iota_iU_TP_i(v\otimes a)$ holds for $v\in\ly$ and $a\in\alg$.
Since $U_T$ and $K_T$ are bounded, $K_T=\sum_{i=1}^{\infty}\iota_iU_TP_i$ holds on $\modu$.
\end{proof}

 \subsection{Decomposition of $K_T$}
We derive a decomposition of $K_T$ called the eigenoperator decomposition.
We first derive a fundamental decomposition using a cocycle on $\mcl{Y}$.
Then, we refine the decomposition using generalized Oseledets subspaces.
\subsubsection{Fundamental decomposition using cocycle}
We first define vectors to decompose the operator $K_T$ using a cocycle on $\mcl{Y}$.
\begin{definition}
For $i\in\mathbb{Z}$, we define a linear operator $w_i:\lz\to\lx$ as
\begin{equation*}
(w_iu)(y,z)=\left\{
\begin{array}{ll}
(U_{g(y,\cdot)}U_{g(h(y),\cdot)}\cdots U_{g(h^{i-1}(y),\cdot)}u)(z)&\quad (i>0)\\
I&\quad (i=0)\\
(U_{g(h^{-1}(y),\cdot)}^*U_{g(h^{-2}(y),\cdot)}^*\cdots U_{g(h^{i}(y),\cdot)}^*u)(z)&\quad (i<0).
\end{array}
\right.
\end{equation*}
\end{definition}

We can see that $\modu$ can be also regarded as a left $\alg$-module.
Thus, we can also consider left $\alg$-linear operators on $\modu$.
In the following, we denote the action of a left $\alg$-linear operator $A$ on a vector $w\in\modu$ by $w\cdot M$.
\begin{proposition}\label{prop:w_i}
For $i\in\mathbb{Z}$, we have $w_i\in\modu$.
Moreover, $K_Tw_i=w_{i+1}=w_i\cdot M_i$, where $M_i$ is a left $\alg$-linear multiplication operator on $\modu$ defined as $(w\cdot M_i)(y)=w(y) U_{g(h^{i}(y),\cdot)}$.
\end{proposition}
\begin{proof}
We obtain $w_i\in\modu$ in the same manner as the proof of Proposition~\ref{prop:K_TinM}.
The identities $K_Tw_i=w_{i+1}=w_i\cdot M_i$ follow by the definition of $w_i$.
\end{proof}

The vectors $w_i$ characterize the dynamics within $\mathcal{Z}$, which are specific to skew product dynamical systems and of particular interest to us.

\begin{proposition}
The action of the Koopman operator $U_T$ is decomposed into two parts as
\begin{equation}
U_T^i(v\otimes u)(y,z)=U_h^iv(z) \cdot U_T^{i-1}u\circ g(y,z)\label{eq:U_Tdec}
\end{equation}
for $v\in\ly$, $u\in\lz$, and $i\in\mathbb{Z}$.
\end{proposition}
\begin{proof}
For $i>0$, it follows by the definition of $U_T$.
Regarding the case of $i\le 0$, $U_T^{-1}$ is calculated as follows:
\begin{align*}
\bracket{v_1\otimes u_1,U_T(v_2\otimes u_2)}
&=\int_{z\in\mcl{Z}}\int_{y\in\mcl{Y}}\overline{v_1(y)u_1(z)}v_2(h(y))u_2(g(y,z))\mr{d}\mu(y)\mr{d}\nu(z)\\
&=\int_{z\in\mcl{Z}}\int_{y\in\mcl{Y}}\overline{v_1(h^{-1}(y))u_1(g_{h^{-1}(y)}(z))}v_2(y)u_2(z)\mr{d}\mu(y)\mr{d}\nu(z),
\end{align*}
where for $y\in\mcl{Y}$, the map $g_y:\mcl{Z}\to\mcl{Z}$ is defined as $g_y(z)=g(y,z)$.
Thus, we have $U_T^{-1}(u\otimes v)(y,z)=v(h^{-1}(y))u(g_{h^{-1}(y)}(z))$.
As a result, the equation~\eqref{eq:U_Tdec} is derived also for $i\le 0$.
\end{proof}
We define a submodule $\mcl{W}$ of $\modu$, which is composed of the vectors $w_i$ ($i\in\mathbb{Z}$).
Let 
\begin{equation*}
\mcl{W}_0=\bigg\{\sum_{i\in F}w_ic_i\,\mid\, F\subseteq\mathbb{Z}:\ \mbox{finite set},\ c_i\in\alg\bigg\}
\end{equation*}
and let $\mcl{W}$ be the completion of $\mcl{W}_0$ with respect to the norm in $\modu$.
Note that $\mcl{W}$ is a submodule of $\modu$ and Hilbert $\alg$-module.
Moreover, for $u\in\lz$, let $\tilde{w}_{u,i}\in\lx$ be defined as $\tilde{w}_{u,i}(y,z)=u(g(h^{i-1}(y),\ldots,g(h(y),g(y,z))\ldots))$ for $i>0$, $\tilde{w}_{u,0}(y,z)=u(z)$, and $\tilde{w}_{u,i}(y,z)=u(g(h^{i}(y),\ldots,g(h^{-2}(y),g(h^{-1}(y),z))\ldots))$ for $i<0$.
Let
\begin{equation*}
\tilde{\mcl{W}}_0=\bigg\{\sum_{j=1}^n\sum_{i\in F}c_i\tilde{w}_{u_j,i}\,\mid\, n\in\mathbb{N},\ F\subseteq\mathbb{Z}:\ \mbox{finite set},\ c_i\in\mathbb{C},\ u_j\in\lz\bigg\}
\end{equation*}
and $\tilde{\mcl{W}}$ be the completion of $\tilde{\mcl{W}}_0$ with respect to the norm in $\lx$.

We show the connection of the operator $K_T$ restricted on $\mcl{W}$ with the Koopman operator $U_T$.
\begin{proposition}
With the notation defined in Proposition~\ref{prop:UandK}, we have $K_T\vert_{\mcl{W}}\,\iota_i\vert_{\tilde{\mcl{W}}}=\iota_i U_T\vert_{\tilde{\mcl{W}}}$ for i=1,2,\ldots.
\end{proposition}
\begin{proof}
For $u\in\lz$, $j\in\mathbb{Z}$, and $i>0$, we have
\begin{align*}
(\iota_i \tilde{w}_{u,j})(y)&=u(g(h^{i-1}(y),\ldots,g(h(y),g(y,\cdot))\ldots))\gamma_i'\\
&=U_{g(h^{i-1}(y),\cdot)}\cdots U_{g(h(y),\cdot)}U_{g(h(y),\cdot)}u\gamma_i'
=w_i(y)(u\gamma_i').
\end{align*}
Thus, we obtain $\iota_i \tilde{w}_{u,j}\in\mcl{W}$.
We obtain $\iota_i \tilde{w}_{u,j}\in\mcl{W}$ for $i\le 0$ in the same manner as the case of $i>0$.
Therefore, the range of $\iota_i\vert_{\tilde{W}}$ is contained in $\mcl{W}$.
The equality is deduced by the definitions of $K_T$ and $U_T$.
\end{proof}

We can describe the decomposition proposed in Proposition~\ref{prop:w_i} using operators on Hilbert $C^*$-modules.
Let 
\begin{align*}
&\mcl{C}_0=\{(\ldots,c_{-1},c_0,c_1,\ldots)\,\mid\,c_i\in\alg,\ c_i=0\mbox{ for all but finite }i\in \mathbb{Z}\},\\
&\mcl{C}'_0=\{(\ldots,A_{-1},A_0,A_1,\ldots)\,\mid\,A_i:\ \mbox{left }\alg\mbox{-linear operator on }\mcl{W},\\
&\qquad\qquad\qquad\qquad\qquad\qquad\qquad A_i=0\mbox{ for all but finite }i\in \mathbb{Z}\}.
\end{align*}
We can see $\mcl{C}_0$ and $\mcl{C}'_0$ are right $\alg$-modules.
We define $\alg$-valued semi-inner products in $\mcl{C}_0$ and $\mcl{C}'_0$ as 
\begin{align*}
&\bracket{(\ldots,c_{-1},c_0,c_{1},\ldots),(\ldots,d_{-1},d_0,d_{1},\ldots)}_{\mcl{C}_0}=\sum_{i,j\in \mathbb{Z}}c_i^*\bracket{w_i,w_j}_{\modu}d_j,\\
&\bracket{(\ldots,A_{-1},A_0,A_1,\ldots),(\ldots,B_{-1},B_0,B_1,\ldots)}_{\mcl{C}'_0}=\sum_{i,j\in \mathbb{Z}}\bracket{w_i\cdot A_i,w_j\cdot B_j}_{\modu},
\end{align*}
respectively.

We define an equivalent relation $c\sim d$ by $c-d\in\mcl{N}$ for $c,d\in\mcl{C}_0$, where $\mcl{N}=\{c\in\mcl{C}_0\,\mid\,\bracket{c,c}_{\mcl{C}_0}=0\}$.
There is an $\alg$-valued inner product on $\mcl{C}_0/\sim$ given by $\bracket{c,d}=\bracket{c+\mcl{N},d+\mcl{N}}$.
We denote by $\mcl{C}$ and $\mcl{C}'$ the completions of $\mcl{C}_0/\sim$ and $\mcl{C}'_0/\sim$ with respect to the norms induced by the above inner products.
We abuse the notation and denote by $(\ldots,c_{-1},c_0,c_{1},\ldots)$ the equivalent class of $(\ldots,c_{-1},c_0,c_{1},\ldots)$ with respect to 
Let $W$ be a right $\alg$-linear operator from $\mcl{C}'$ to $\mcl{W}$ defined as
\begin{equation*}
W(\ldots,A_{-1},A_0,A_{1},\ldots)=\sum_{i\in \mathbb{Z}}w_i\cdot A_i    
\end{equation*}
for $(\ldots,A_{-1},A_0,A_{1},\ldots)\in\mcl{C}'_0$ and let $X$ be a right $\alg$-linear operator from $\mcl{W}$ to $\mcl{C}$ defined as 
\begin{equation*}
X\sum_{i\in F}w_ic_i=(\ldots,c_{-1},c_0,c_{1},\ldots)
\end{equation*}
for a finite set $F\subseteq\mathbb{Z}$.
In addition, let $M$ be a right $\alg$-linear operator from $\mcl{C}$ to $\mcl{C}'$ defined as
\begin{equation*}
M(\ldots,c_{-1},c_0,c_{1},\ldots)=(\ldots,M_{-1}c_{-1},M_0c_{0},M_{1}c_1,\ldots)
\end{equation*}
for $(\ldots,c_{-1},c_0,c_{1},\ldots)\in\mcl{C}_0$, which is formally denoted by $\opn{diag}\{M_i\}_{i\in\mathbb{Z}}$.
Here, $M_i$ is the multiplication operator defined in Proposition~\ref{prop:w_i}.
\begin{proposition}
The operators $W$ and $X$ are unitary operators.
Therefore, $XW$ is an unitary operator from $\mcl{C}'$ to $\mcl{C}$.
\end{proposition}
\begin{proof}
Let $\tilde{W}:\mcl{W}\to\mcl{C}'$ be the right $\alg$-linear operator defined as $\tilde{W}(\sum_{i\in F}w_ic_i)=(\ldots,C_{-1},C_0,C_1,\ldots)$, where $C_i$ is the left $\alg$-linear multiplication operator on $\mcl{W}$ with respect to the constant function $c_i$.
In addition, let $\tilde{X}:\mcl{C}\to\mcl{W}$ be the right $\alg$-linear operator defined as $\tilde{X}(\ldots,c_{-1},c_0,c_1,\ldots)=\sum_{i\in F}w_ic_i$.
Then, $\tilde{W}$ and $\tilde{X}$ are the inverses of $W$ and $X$, respectively.
Moreover, for $(\ldots,A_{-1},A_0,A_{1},\ldots), (\ldots,B_{-1},B_0,B_{1},\ldots)\in\mcl{C}'_0$, we have
\begin{align*}
&\bracket{W(\ldots,A_{-1},A_0,A_{1},\ldots),W(\ldots,B_{-1},B_0,B_{1},\ldots)}_{\modu}
=\Bbracket{\sum_{i\in\mathbb{Z}}w_i\cdot A_i,\sum_{i\in\mathbb{Z}}w_i\cdot B_i}_{\modu}\\
&\qquad =\bracket{(\ldots,A_{-1},A_0,A_{1},\ldots),(\ldots,B_{-1},B_0,B_{1},\ldots)}_{\mcl{C}'}
\end{align*}
and for $c_i,d_i\in\alg$, finite subsets $F$ and $G$ of $\mathbb{Z}$, we have
\begin{align*}
\Bbracket{X\sum_{i\in F}w_ic_i,X\sum_{i\in G}w_id_i}_{\mcl{C}}
&=\bracket{(\ldots,c_{-1},c_0,c_{1},\ldots),(\ldots,d_{-1},d_0,d_{1},\ldots)}_{\mcl{C}}\\
&=\Bbracket{\sum_{i\in F}w_ic_i,\sum_{i\in G}w_id_i}_{\modu}.
\end{align*}
\end{proof}
\begin{proposition}
The operator $M$ is well-defined and $K_T\vert_{\mcl{W}}=WMX$.
\end{proposition}
\begin{proof}
Since $K_T$ is unitary, we have
\begin{equation}
\bracket{w_{i+1},w_{j+1}}=\bracket{K_Tw_i,K_Tw_j}=\bracket{w_i,w_j}\label{eq:w_i_w_i+1}
\end{equation}
for $i,j\in\mathbb{Z}$.
Assume $(\ldots,c_{-1},c_0,c_1,\ldots)=0$.
Then, by the equation~\eqref{eq:w_i_w_i+1}, we obtain
\begin{align*}
&\bracket{M(\ldots,c_{-1},c_0,c_1,\ldots),M(\ldots,c_{-1},c_0,c_1,\ldots)}_{\mcl{C}'}\\
&\qquad=\bracket{(\ldots,M_{-1}c_{-1},M_0c_0,M_1c_1,\ldots),(\ldots,M_{-1}c_{-1},M_0c_0,M_1c_1,\ldots)}_{\mcl{C}'}\\
&\qquad=\sum_{i,j\in F}\bracket{w_i\cdot M_ic_i,w_j\cdot M_jc_j}_{\modu}
=\sum_{i,j\in F}\bracket{w_{i+1}c_i,w_{j+1}c_j}_{\modu}\\
&\qquad=\sum_{i,j\in F}\bracket{w_{i}c_i,w_{j}c_j}_{\modu}
=0,
\end{align*}
which shows the well-definiteness of $M$.
The decomposition $K_T\vert_{\mcl{W}}=WMX$ is derived by Proposition~\ref{prop:w_i}.
\end{proof}
In summary, we obtain the following commutative diagram:
\begin{equation*}
 \xymatrix{
\mcl{W}\ar[rr]^{K_T}\ar[d]<0.5ex>^{X}&  &\mcl{W}\ar[d]<0.5ex>^{W^*}\\
\mcl{C} \ar[rr]^{M}\ar[u]<0.5ex>^{X^*}& &\mcl{C}' \ar[u]<0.5ex>^{W}
}
\end{equation*}

\subsubsection{Further decomposition}
We further decompose $w_i$ and $M_i$ and obtain a more detailed decomposition of $K_T\vert_{\mcl{W}}$.
Let $\mcl{V}_1,\mcl{V}_2,\ldots$ be a sequence of maps from $\mcl{Y}$ to the set of all closed subspaces of $\lz$  satisfying $\lz=\overline{\opn{Span}\{\bigcup_{j=1}^{\infty}\mcl{V}_j(y)\}}$ for a.s. $y\in\mcl{Y}$.
Let $p_{j}(y)\in\alg$ be the projection onto $\mcl{V}_j(y)$, i.e., it satisfies $p_j(y)^2=p_j(y)$ and $p_j(y)^*=p_j(y)$.
For $i\in\mathbb{Z}$ and $j=1,2,\ldots$, we define a linear map $\hat{w}_{i,j}$ from $\lz$ to $\lx$ as $(\hat{w}_{i,j}u)(y,z)=(w_i(y)p_{j}(h^i(y))u)(z)$.
We decompose $K_T\vert_{\mcl{W}}$ using $\hat{w}_{i,j}$.
For each $j=1,2,\ldots$, the following theorem holds:
\begin{theorem}[Eigenoperator decomposition for discrete-time systems]\label{prop:detailed_decom}
Assume $\mcl{V}_j$ satisfies $U_{g(y,\cdot)}\mcl{V}_j(h(y))\subseteq\mcl{V}_j(y)$ for a.s. $y\in\mcl{Y}$.
Assume in addition, for any $u\in\lz$ and $i\in\mathbb{Z}$, the map $(y,z)\mapsto (p_{j}(h^i(y))u)(z)$ is measurable.
Then, $\hat{w}_{i,j}$ is contained in $\modu$ and we have $K_T\hat{w}_{i,j}=\hat{w}_{i+1,j}=\hat{w}_{i,j}\cdot \hat{M}_{i,j}$.
Here, $\hat{M}_{i,j}$ is a left $\alg$-linear multiplication operator on $\modu$ defined as $(w\cdot \hat{M}_{i,j})(y)=w(y) U_{g(h^{i}(y),\cdot)}p_j(h^i(y))$.
\end{theorem}

\begin{proof}
We have
\begin{equation*}
K_T\hat{w}_{i,j}(y)=U_{g(y,\cdot)}\hat{w}_{i,j}(h(y))=w_{i+1}(y)p_{j}(h^{i+1}(y))=\hat{w}_{i+1,j}(y).
\end{equation*}
In addition, since by the assumption, the range of $U_{g(h^i(y),\cdot)}p_{j}(h^{i+1}(y))$ is contained in $\mcl{V}(h^i(y))$, we have
\begin{align*}
K_T\hat{w}_{i,j}(y)&=w_{i}(y)U_{g(h^i(y),\cdot)}p_{j}(h^{i+1}(y))\\
&=w_{i}(y)p_j(h^i(y))U_{g(h^i(y),\cdot)}p_{j}(h^{i+1}(y))=\hat{w}_{i,j}(y)\cdot\hat{M}_{i,j}.
\end{align*}
\end{proof}
\begin{corollary}\label{cor:detailed_dec}
By replacing $\{w_i\}_{i\in\mathbb{Z}}$ and $\{M_i\}_{i\in\mathbb{Z}}$ by $\{\hat{w}_{i,j}\}_{i\in\mathbb{Z},j\in\mathbb{N}}$ and $\{\hat{M}_{i,j}\}_{i\in\mathbb{Z},j\in\mathbb{N}}$, respectively,
we define $\hat{\mcl{W}}$, $\hat{\mcl{C}}$, $\hat{\mcl{C}'}$, $\hat{W}$, $\hat{X}$, and $\hat{M}$ in the same manner as $\mcl{W}$, $\mcl{C}$, $\mcl{C}'$, $W$, $X$, and $M$, respectively.
Then, under the assumptions of Theorem~\ref{prop:detailed_decom}, we obtain $K_T\vert_{\mcl{W}}=\hat{W}\hat{M}\hat{X}$.
\end{corollary}

We call $\hat{M}_{i,j}$ an eigenoperator and $\hat{w}_{i,j}$ an eigenvector.
In addition, we call the subspace $\mcl{V}_j(y)$ satisfying the assumption in Theorem~\ref{prop:detailed_decom} generalized Oseledets space.

If $\mcl{V}_j(y)$ is a finite-dimensional space, then we can explicitly calculate the spectrum of $\hat{M}_{i,j}$ as follows.
\begin{proposition}\label{prop:spectrum_finitedim}
Assume $\opn{dim}(\mcl{V}_j(y))$ is finite and constant with respect to $y \in \mcl{Y}$.
Then, $\sigma(\hat{M}_{i,j})=\{\lambda\in\mathbb{C}\,\mid\,^{\forall}\epsilon>0,\ \mu(\{y\in\mcl{Y}\,\mid\,\lambda\in\sigma_{\epsilon}(f_{i,j}(y))\})>0\}$.
Here, $f_{i,j}(y)=U_{g(h^{i}(y),\cdot)}p_j(h^i(y))$. 
In addition, $\sigma(a)$ and $\sigma_{\epsilon}(a)$ for $a\in\alg$ are the spectrum and the essential spectrum of $a$, respectively.
\end{proposition}
\begin{proof}
For $\lambda\in\mathbb{C}$, $\lambda I-\hat{M}_{i,j}$ is not invertible if and only if there exists $w\in\modu$ with $w\neq 0$ such that $w(y)(\lambda I-f_{i,j}(y))=0$ for a.s. $y\in\mcl{Y}$, which is equivalent to
\begin{align*}
\mu(\{y\in\mcl{Y}\,\mid\,\lambda I-f_{i,j}(y)\mbox{ is not invertible}\})>0.
\end{align*}

Assume $\lambda I-f_{i,j}(y)$ is invertible for a.s. $y\in\mcl{Y}$.
Then, we have $(\lambda I-\hat{M}_{i,j})^{-1}w(y)=w(y)(\lambda I-f_{i,j}(y))^{-1}$.
For $w\in\modu$, we have
\begin{align*}
\Vert (\lambda I-\hat{M}_{i,j})^{-1}w\Vert_{\modu}^2
&=\bigg\Vert \int_{\mcl{Y}}(\lambda I-f_{i,j}(y))^{-*}w(y)^{*}w(y)(\lambda I-f_{i,j}(y))^{-1}\mr{d}\mu(y)\bigg\Vert_{\alg}\\
&\le \opn{tr}\bigg(\int_{\mcl{Y}}(\lambda I-f_{i,j}(y))^{-*}w(y)^{*}w(y)(\lambda I-f_{i,j}(y))^{-1}\mr{d}\mu(y)\bigg)\\
&=\opn{tr}\bigg(\int_{\mcl{Y}}w(y)(\lambda I-f_{i,j}(y))^{-*}(\lambda I-f_{i,j}(y))^{-1}w(y)^{*}\mr{d}\mu(y)\bigg).
\end{align*}
Thus, we have
\begin{align*}
&\frac{1}{d_j}\bigg\Vert \int_{\mcl{Y}}w(y)(\lambda I-f_{i,j}(y))^{-*}(\lambda I-f_{i,j}(y))^{-1}w(y)^{*}\mr{d}\mu(y)\bigg\Vert_{\alg}
 \le\Vert (\lambda I-\hat{M}_{i,j})^{-1}w\Vert_{\modu}^2\\
&\qquad \le d_j\bigg\Vert \int_{\mcl{Y}}w(y)(\lambda I-f_{i,j}(y))^{-*}(\lambda I-f_{i,j}(y))^{-1}w(y)^{*}\mr{d}\mu(y)\bigg\Vert_{\alg},
\end{align*}
where $d_j=\opn{dim}(\mcl{V}_j(y))$.
Assume for any $\epsilon>0$, $\mu(\{y\in\mcl{Y}\,\mid\,\Vert (\lambda I-f_{i,j}(y))^{-1}\Vert_{\alg}>1/\epsilon\})>0$.
We set $a_{i,j}(y)=v_{i,j}(y)v_{i,j}(y)^*$, where $v_{i,j}(y)$ is the orthonormal eigenvector corresponding to the largest eigenvalue of $(\lambda I-f_{i,j}(y))^{-*}\lambda I-f_{i,j}(y))^{-1}$.
Then, we have
\begin{align*}
&\Vert (\lambda I-\hat{M}_{i,j})^{-1}w\Vert_{\modu}^2
\ge \frac{1}{d_j}\bigg\Vert \int_{\mcl{Y}}\Vert(\lambda I-f_{i,j}(y))^{-1}\Vert_{\alg}^2 w(y)a_{i,j}(y)w(y)^{*}\mr{d}\mu(y)\bigg\Vert_{\alg}\\
&\qquad\ge\mu(\{y\in\mcl{Y}\,\mid\,\Vert (\lambda I-f_{i,j}(y))^{-1}\Vert_{\alg}>1/\epsilon\})\frac{1}{\epsilon^2 d_j}\bigg\Vert \int_{\mcl{Y}} w(y)a(y)w(y)^{*}\mr{d}\mu(y)\bigg\Vert_{\alg}.
\end{align*}
Setting $w(y)=a_{i,j}(y)$, we derive that $(\lambda I-\hat{M}_{i,j})^{-1}$ is unbounded.
Conversely, assume $(\lambda I-\hat{M}_{i,j})^{-1}$ is unbounded.
Then, we obtain
\begin{align*}
&\Vert (\lambda I-\hat{M}_{i,j})^{-1}w\Vert_{\modu}^2
\le d_j\esssup_{y\in\mcl{Y}}\Vert (\lambda I-f_{i,j}(y))^{-1}\Vert^2 \bigg\Vert\int_{\mcl{Y}}w(y)^*w(y)\mr{d}\mu(y)\bigg\Vert_{\alg}\\
&\qquad\le d_j\esssup_{y\in\mcl{Y}}\Vert (\lambda I-f_{i,j}(y))^{-1}\Vert^2 \Vert w\Vert_{\modu}^2.
\end{align*}
Thus, for any $\epsilon>0$, $\mu(\{y\in\mcl{Y}\,\mid\,\Vert (\lambda I-f_{i,j}(y))^{-1}\Vert_{\alg}>1/\epsilon\})>0$, which completes the proof.
\end{proof}

\subsubsection{Construction of the generalized Oseledets space $\mcl{V}_{j}(y)$}
For cocycles generated by matrices, the existence of the Oseledets space is guaranteed by the multiplicative ergodic theorem.
This theorem has been generalized for cocycles generated by compact operators or operators that have similar properties to the compactness~\cite{thieullen87,Schaumloffel91}.
In our case, since the cocycle is generated by a unitary operator, we can construct $\mcl{V}_j$ explicitly if $h$ is periodic.
\begin{proposition}\label{prop:construct_osedelets}
Assume $h^n(y)=y$ for any $y\in\mcl{Y}$.
Let $U:\mcl{Y}\to\mcl{B}(\oplus_{k=1}^n\lz)$ be defined as 
\begin{equation*}
U(y)=(U_{g(h(y),\cdot)}\oplus U_{g(h^2(y),\cdot)}\oplus\cdots\oplus U_{g(h^{n-1}(y),\cdot)}\oplus U_{g(y,\cdot)})S, 
\end{equation*}
where $S$ is the permutation operator defined as $S\oplus_{k=1}^nu_k=\oplus_{k=1}^{n-1}u_{k+1}\oplus u_1$.
Let $E(y)$ be the spectral measure with respect to $U(y)$ and $\mcl{T}\subset [0,2\pi)$ be a subset of $[0,2\pi)$.
Let $\tilde{\mcl{V}}_{k}(y)$ be the range of $P_kE(y)(T)$ and let $\mcl{V}(y)=\tilde{\mcl{V}}_{n}(y)$, where $P_k:\oplus_{l=1}^n\lz\to\lz$ is the projection defined as $\oplus_{l=1}^nu_l\mapsto u_k$.
Then, we have $U_{g(y,\cdot)}\mcl{V}(h(y))\subseteq\mcl{V}(y)$.
\end{proposition}
\begin{proof}
We first show $\oplus_{k=1}^n \tilde{\mcl{V}}_{k}(y)$ is an invariant subspace of $U(y)$.
Let $u\in E(y)(\mcl{T})$. Then, we have
\begin{equation*}
U(y)\tilde{P}_ku=U_{g(h^{k-1}(y),\cdot)}u_k=\tilde{P}_kU(y)u
\end{equation*}
for $k=1,\ldots,n$.
Here, $\tilde{P}_k:\oplus_{l=1}^n\lz\to\oplus_{l=1}^n\lz$ is the linear operator defined as $\oplus_{l=1}^nu_l\mapsto \oplus_{l=1}^{n}\tilde{u}_l$, where $\tilde{u}_l=0$ for $l\neq k$ and $\tilde{u}_l=u_k$ for $l=k$.
Since $E(y)(\mcl{T})$ is an invariant subspace of $U(y)$, we have $U(y)u\in E(y)(\mcl{T})$.
Thus, we have $U(y)\oplus_{k=1}^n \tilde{\mcl{V}}_{k}(y)\subseteq \oplus_{k=1}^n \tilde{\mcl{V}}_{k}(y)$.
Therefore, we have $U_{g(y,\cdot)}\tilde{\mcl{V}}_{n}(h(y))\subseteq \tilde{\mcl{V}}_{n-1}(h(y))$.
In addition, for $\oplus_{k=1}^nu_k\in \oplus_{k=1}^n\lz$, we have 
\begin{align*}
S^{-1}U(h(y))S\oplus_{k=1}^nu_k&=S^{-1}U(h(y))\oplus_{k=1}^{n}u_{k+1}
=S^{-1}\oplus_{k=1}^{n}U_{g(h^{k+1}(y),\cdot)}(y)u_{k+2}\\
&=\oplus_{k=1}^{n}U_{g(h^{k}(y),\cdot)}(y)u_{k+1}
=U(y)\oplus_{k=1}^nu_k,
\end{align*}
where $u_{k+n}=u_k$ for $k=1,2$.
Thus, we have $S^{-1}U(h(y))S=U(y)$, and the spectral measure $E(h(y))$ of $U(h(y))$ is represented as $E(h(y))=SE(y)S^{-1}$.
Therefore, for $k=1,2,\ldots$, we obtain
\begin{equation*}
P_kE(h(y))(\mcl{T})=P_kSE(y)(\mcl{T})S^{-1}=P_{k+1}E(y)(\mcl{T})S^{-1},
\end{equation*}
where $P_{k+1}=P_1$, which implies $\tilde{\mcl{V}}_{k}(h(y))=\tilde{\mcl{V}}_{k+1}(y)$.
Combining this identity with the inclusion $U_{g(y,\cdot)}\mcl{V}_{n}(h(y))\subseteq \mcl{V}_{n-1}(h(y))$, we have $U_{g(y,\cdot)}\tilde{\mcl{V}}_{n}(h(y))\subseteq \tilde{\mcl{V}}_{n}(y)$.
\end{proof}

\begin{corollary}
Let $\mcl{T}_1,\mcl{T}_2,\ldots\subset [0,2\pi)$ be a sequence of countable disjoint subsets of $[0,2\pi)$ such that $[0,2\pi)=\bigcup_{j=1}^{\infty}\mcl{T}_j$.
If we set $\mcl{V}_j(y)$ as $\mcl{V}(y)$ in Proposition~\ref{prop:construct_osedelets} by replacing $\mcl{T}$ with $\mcl{T}_j$, it satisfies the assumption in Theorem~\ref{prop:detailed_decom}.
\end{corollary}

\subsubsection{Connection with Koopman operator on Hilbert space}
Assume $p_j(y)=\hat{p}_j$ for any $y\in\mcl{Y}$, where $\hat{p}_j\in\alg$ is a projection.
By Proposition~\ref{prop:detailed_decom}, we obtain the following commutative diagram:
\begin{equation*}
 \xymatrix{
\mcl{M}\ar[rr]^{K_T}\ar[d]^{P_j}&  &\mcl{M}\ar[d]^{P_j}\\
\ly\otimes\hat{\alg}_j \ar[rr]^{K_T}& &\ly\otimes\hat{\alg}_j 
}
\end{equation*}
where $\hat{\alg}_j=\{\hat{p}_ja\hat{p}_j\,\mid\,a\in\alg\}$ is a $C^*$-subalgebra of $\alg$ and $P_j:\modu\to\ly\otimes\hat{\alg}_j$ is defined as $w\cdot P_j=w\hat{p}_j$.
If $\mcl{V}_j$ is a finite dimensional space, then $\ly\otimes\hat{\alg}_j$ is isomorphic to a Hilbert space and the action of $K_T$ on $\ly\otimes\hat{\alg}_j$ is reduced to that of $U_T$.
\begin{proposition}\label{prop:finite_dim}
Assume $\mcl{V}_j$ is an $n$-dimensional space.
Let $\{\gamma_1,\ldots,\gamma_n\}$ be an orthonormal basis of $\mcl{V}_j$ and let $\lambda_j:\ly\otimes\hat{\alg}_j\to\oplus_{i=1}^n\lx$ be a linear operator defined as $\lambda_j(v\otimes a)=(v\otimes a\gamma_1,\ldots,v\otimes a\gamma_n)$.
Then, $\lambda_j$ is an isomorphism and we have the following commutative diagram: 
\begin{equation*}
 \xymatrix{
\ly\otimes\hat{\alg}_j\ar[rr]^{K_T}\ar[d]^{\lambda_j}&  &\ly\otimes\hat{\alg}_j\ar[d]^{\lambda_j}\\
 \oplus_{i=1}^n\lx\ar[rr]^{\oplus_{i=1}^nU_T}& &\oplus_{i=1}^n\lx 
}
\end{equation*}
\end{proposition}
\begin{proof}
Let $\tilde{\lambda}_j:\oplus_{i=1}^n\lx\to\ly\otimes\hat{\alg}_j$ be a linear operator defined as $\tilde{\lambda}_j(v_1\otimes u_1,\ldots,v_n\otimes u_n)=\sum_{i=1}^nv_i\otimes u_i\gamma_i'$.
Then, $\tilde{\lambda}_j$ is the inverse of $\lambda_j$.
In addition, we have
\begin{align*}
&\bracket{\lambda_j(v\otimes a),\lambda_j(v\otimes a)}_{\oplus_{i=1}^n\lx}
=\sum_{i=1}^n\bracket{v\otimes a\gamma_i,v\otimes a\gamma_i}_{\lx}\\
&\qquad=\bracket{v,v}_{\ly}\sum_{i=1}^n\bracket{\gamma_i,a^*a\gamma_i}_{\lz}
\le n\Vert\bracket{v,v}_{\ly}a^*a\Vert_{\alg}
=n\Vert v\otimes a\Vert_{\modu}^2.
\end{align*}
Thus, $\lambda_j$ is an isomorphism.
The commutativity of the diagram is derived by Proposition~\ref{prop:UandK}.
\end{proof}

\section{Examples}
\label{sec:discrete_time_examples}

\subsection{The case of $\mcl{Z}$ is a compact Hausdorff group}
Let $\mcl{Z}$ be a compact Hausdorff group equipped with the (normalized) Haar measure $\nu$.
Let $\hat{\mcl{Z}}$ be the set of equivalent classes of irreducible unitary representations. 
For an irreducible representation $\rho$, let $\mcl{E}_{\rho}$ be the representation space of $\rho$ and let $n_{\rho}$ be the dimension of $\mcl{E}_{\rho}$.
Note that since $\mcl{Z}$ is a compact group, $n_{\rho}$ is finite.
Let $\{e_{\rho,1},\ldots,e_{\rho,n_{\rho}}\}$ be an orthonormal basis of $\mcl{E}_{\rho}$ and let $\gamma_{\rho,i,j}:\mcl{Z}\to\mathbb{C}$ be the matrix coefficient defined as $\gamma_{\rho,i,j}(z)=\bracket{e_{\rho,i},\rho(z)e_{\rho,j}}$.
By the Peter--Weyl theorem, $\bigcup_{[\rho]\in\zhat}\{\gamma_{\rho,i,j}\;\mid\;i,j=1,\ldots,n_{\rho}\}$ is an orthonormal basis of $\lz$, where $[\rho]$ is the equivalent class of an irreducible representation $\rho$.
We set the map $g:\mcl{Y}\times\mcl{Z}\to\mcl{Z}$ as $g(y,z)=z\tilde{g}(y)$, where $\tilde{g}:\mcl{Y}\to\mcl{Z}$ is a measurable map.
Let $\Gamma_{\rho,i}:\mcl{E}_{\rho}\to\lz$ be the linear operator defined as $e_{\rho,j}\mapsto \gamma_{\rho,i,j}$ for $i,j=1,\ldots,n_{\rho}$.
Note that the adjoint $\Gamma_{\rho,i}^*:\lz\to\mcl{E}_{\rho}$ is written as $u\mapsto \sum_{j=1}^{n_\rho}\bracket{\gamma_{\rho,i,j},u}e_{\rho,j}$.
Then, regarding the Koopman operator $U_{g(y,\cdot)}$ on $\lz$, we have
\begin{align*}
&U_{g(y,\cdot)}\Gamma_{\rho,i}\Gamma_{\rho,i}^*u(z)
=U_{g(y,\cdot)}\sum_{j=1}^{n_\rho}\bracket{\gamma_{\rho,i,j},u}\gamma_{\rho,i,j}(z)
=\sum_{j=1}^{n_\rho}\bracket{\gamma_{\rho,i,j},u}\gamma_{\rho,i,j}(\tilde{g}(y)z)\\
&\qquad=\sum_{j=1}^{n_\rho}\bracket{\gamma_{\rho,i,j},u}\bracket{e_{\rho,i},\rho(z\tilde{g}(y))e_{\rho,j}}
=\sum_{j=1}^{n_\rho}\bracket{\gamma_{\rho,i,j},u}\bracket{\rho(z)^*e_{\rho,i},\rho(\tilde{g}(y))e_{\rho,j}}\\
&\qquad=\sum_{j=1}^{n_\rho}\bracket{\gamma_{\rho,i,j},u}\bracket{\rho(z)^*e_{\rho,i},\sum_{k=1}^{n_{\rho}}\bracket{e_{\rho,k},\rho(\tilde{g}(y))e_{\rho,j}}e_{\rho,k}}\\
&\qquad=\sum_{j,k=1}^{n_\rho}\bracket{\gamma_{\rho,i,j},u}\bracket{\rho(\tilde{g}(y))^*e_{\rho,k},e_{\rho,j}}\gamma_{\rho,i,k}(z)\\
&\qquad=\sum_{k=1}^{n_\rho}\bracket{\rho(\tilde{g}(y))^*e_{\rho,k},\Gamma_{\rho,i}^*u}\gamma_{\rho,i,k}(z)
=\bigg(\Gamma_{\rho,i}\sum_{k=1}^{n_\rho}\bracket{\rho(\tilde{g}(y))^*e_{\rho,k},\Gamma_{\rho,i}^*u}e_{\rho,k}\bigg)(z)\\
&\qquad=\Gamma_{\rho,i}\rho(\tilde{g}(y))\Gamma_{\rho,i}^*u(z)
\end{align*}
for $u\in\lz$, $z\in\mcl{Z}$, and $i=1,\ldots,n_{\rho}$.
Thus, we have $U_{g(y,\cdot)}=\sum_{[\rho]\in\zhat}\sum_{i=1}^{n_\rho}\Gamma_{\rho,i}\rho(\tilde{g}(y))\Gamma_{\rho,i}^*$.
Therefore, the range of $\Gamma_{\rho,i}$ is an invariant subspace of $U_{g(y,\cdot)}$ for any $y\in\mcl{Y}$.
Thus, we set $\mcl{V}_{[\rho],j}$ as the constant map which takes its value the range of $\Gamma_{\rho,j}$, and apply Proposition~\ref{prop:detailed_decom}.
In this case, the multiplication operator $\hat{M}_{i,[\rho],j}$ is calculated as $(w\cdot\hat{M}_{i,[\rho],j})(y)=w(y)\Gamma_{\rho,j}\rho(\tilde{g}(h^i(y)))\Gamma_{\rho,j}^*$, and by Proposition~\ref{prop:spectrum_finitedim}, its spectrum is calculated as 
\begin{align*}
\sigma(\hat{M}_{i,[\rho],j})&=\{\lambda\in\mathbb{C}\;\mid\;^{\forall}\epsilon>0,\  \mu(\{y\in\mcl{Y}\;\mid\;\lambda\in \sigma_{\epsilon}(\Gamma_{\rho,j}\rho(\tilde{g}(h^i(y)))\Gamma_{\rho,j}^*)\})>0\}
\end{align*}
Note that since $\rho(\tilde{g}(h^i(y)))$ is a linear operator on a finite dimensional space, it has only point spectra.
By Corollary~\ref{cor:detailed_dec}, we obtain a discrete decomposition of $K_T\vert_{\mcl{W}}$ with the multiplication operators $\hat{M}_{i,[\rho],j}$.
Let $\hat{p}_{[\rho],j}=\Gamma_{\rho,j}\Gamma_{\rho,j}^*$.
Then, $K_T$ maps $\ly\otimes\hat{\alg}_{[\rho],j}$ to $\ly\otimes\hat{\alg}_{[\rho],j}$, where $\hat{\alg}_{[\rho],j}=\{\hat{p}_{[\rho],j}a\hat{p}_{[\rho],j}\;\mid\;a\in\alg\}$.
Since $\mcl{V}_{[\rho],j}$ is a finite dimensional space, by Proposition~\ref{prop:finite_dim}, the action of $K_T$ restricted to $\ly\otimes\hat{\alg}_{[\rho],j}$ is reduced to that of $\otimes_{i=1}^{n_{\rho}}U_T$ on $\otimes_{i=1}^{n_{\rho}}(\lx)$ as $K_T=\lambda_{j,[\rho]}(\otimes_{i=1}^{n_{\rho}}U_T)\lambda_{[\rho],j}^{-1}$, where $\lambda_{[\rho],j}(v\otimes a)=(v\otimes a\gamma_{\rho,j,1},\ldots,v\otimes a\gamma_{\rho,j,n_{\rho}})$.

\subsection{The case of $\mcl{Z}=\mathbb{Z}$}
Let $\mcl{Z}=\mathbb{Z}$ equipped with the counting measure.
We set the map $g:\mcl{Y}\times \mcl{Z}\to\mcl{Z}$ as $g(y,z)=z+\tilde{g}(y)$, where $\tilde{g}:\mcl{Y}\to\mcl{Z}$ is a measurable map.
For $i\in\mathbb{Z}$, let $e_i:\mathbb{T}\to\mathbb{C}$ be defined as $e_i(\omega)=\mr{e}^{\sqrt{-1}i\omega}$, where $\mathbb{T}=\mathbb{R}/2\pi\mathbb{Z}$.
Note that $\{e_i\;\mid\;i\in\mathbb{Z}\}$ is an orthonormal basis of $L^2(\mathbb{T})$.
In addition, for $i\in\mathbb{Z}$, let $\gamma_i:\mathbb{Z}\to\mathbb{C}$ be defined as $\gamma_i(i)=1$, $\gamma_i(z)=0\ (z\neq i)$.
Note also that $\{\gamma_i\;\mid\;i\in\mathbb{Z}\}$ is an orthonormal basis of $L^2(\mathbb{Z})$.
Let $\Gamma:L^2(\mathbb{T})\to L^2(\mathbb{Z})$ be the linear operator defined as $e_i\mapsto\gamma_i$ for any $i\in\mathbb{Z}$ and let $\phi_y(\omega)=\mr{e}^{\sqrt{-1}\tilde{g}(y)\omega}$.
Then, we have
\begin{align*}
\Gamma M_{\phi_y}\Gamma^*\gamma_i=\Gamma M_{\phi_y}e_i=\Gamma e_i\phi
=\Gamma e_{i+\tilde{g}(y)}=\gamma_{i+\tilde{g}(y)}=U_{g(y,\cdot)}\gamma_i,
\end{align*}
where $M_{\phi_y}$ is the multiplication operator on $L^2(\mathbb{T})$ defined as $M_{\phi_y}u(\omega)=u(\omega)\phi_y(\omega)=u(\omega)\mr{e}^{\sqrt{-1}\tilde{g}(y)\omega}$.
Thus, we have the spectral decomposition $U_{g(y,\cdot)}=\int_{\omega\in\mathbb{T}}\mr{e}^{\sqrt{-1}\tilde{g}(y)\omega}\mr{d}E(\omega)$, where $E$ is the spectral measure defined as $E(\Omega)=\Gamma M_{\chi_{\Omega}}\Gamma^*$ for a Borel set $\Omega$ and $\chi_{\Omega}$ is the characteristic function of $\Omega$.
Let $T_1,T_2,\ldots$ be a sequence of countable disjoint subsets of $\mathbb{T}$ such that $\mathbb{T}=\bigcup_{j=1}^{\infty}T_j$.
Then, the range of $E(T_j)$ is an invariant subspace of $U_{g(y,\cdot)}$ for any $y\in\mcl{Y}$.
Thus, we set $\mcl{V}_j$ as the constant map which takes its value the range of $E(T_j)$, and apply Proposition~\ref{prop:detailed_decom}.
In this case, $\hat{M}_{i,j}$ is calculated as $(w\cdot\hat{M}_{i,j})(y)=w(y)\Gamma M_{\phi_{h^i(y)}}M_{\chi_{T_j}}\Gamma^*$. 
Let $\hat{p}_{j}=E(T_j)$.
Then, $K_T$ maps $\ly\otimes\hat{\alg}_{j}$ to $\ly\otimes\hat{\alg}_{j}$. 
Since $\mcl{V}_j$ is an infinite dimensional space, we cannot reduce the action of $K_T$ restricted to $\ly\otimes\hat{A}_j$ to that of $U_T$ on a Hilbert space.
However, by Corollary~\ref{cor:detailed_dec}, we obtain a discrete decomposition of $K_T\vert_{\mcl{W}}$ in the Hilbert $C^*$-module even in this case of the spectral decomposition of $U_{g(y,\cdot)}$ is continuous.

\section{Continuous-time systems}
\label{sec:continuous_time_systems}
\subsection{Skew product system and Koopman operator on Hilbert space}
As in the Section~\ref{sec:discrete_systems}, let $\mcl{Y}$ and $\mcl{Z}$ be separable measure spaces equipped with measures $\mu$ and $\nu$, respectively
and let $\mcl{X}=\mcl{Y}\times\mcl{Z}$, the direct product measure space of $\mcl{Y}$ and $\mcl{Z}$.
Let $h:\mathbb{R}\times\mcl{Y}\to\mcl{Y}$ be a map such that for any $t\in\mathbb{R}$, $h(t,\cdot)$ is a measure preserving and invertible map on $\mcl{Y}$.
Moreover, let $g:\mathbb{R}\times \mcl{X}\to\mcl{Z}$ be a map such that for any $t\in\mathbb{R}$, $g(t,\cdot,\cdot)$ is a measurable map from $\mcl{X}$ to $\mcl{Z}$ and for any $y\in\mcl{Y}$, $g(t, y,\cdot)$ is measure preserving and invertible on $\mcl{Z}$.
Consider the following skew product flow on $\mcl{X}$:
\begin{equation*}
 \Phi(t,y,z)=(h(t,y),g(t,y,z))
\end{equation*}
that satisfies $\Phi(0,y,z)=(y,z)$ and $\Phi(s,\Phi(t,y,z))=\Phi(s+t,y,z)$ for any $s,t\in\mathbb{R}$, $y\in\mcl{Y}$, and $z\in\mcl{Z}$.
We denote $\Phi(t,\cdot,\cdot)=\Phi_t$, $h(t,\cdot)=h_t$, and $g(t,\cdot,\cdot)=g_t$, respectively.
For $t\in\mathbb{R}$, we consider the Koopman operator $U_{\Phi_t}$ on $L^2(\mcl{X})$.
Instead of $U_T$ for discrete systems, we consider a family of Koopman operators $\{U_{\Phi_t}\}_{t\in\mathbb{R}}$ for continuous systems.

\subsection{Operator on Hilbert $C^*$-module related to the Koopman operator}
Analogous to the case of discrete systems, we extend the Koopman operator $U_{\Phi_t}$ to an operator on the Hilbert $C^*$-module $\modu$.
\begin{definition}
For $t\in\mathbb{R}$, we define a right $\alg$-linear operator $K_{\Phi_t}$ on $\modu$ by
\begin{equation*}
K_{\Phi_t}(v\otimes a)(y)=v(h_t(y))U_{g_t(y,\cdot)}a
\end{equation*}
for $v\in\ly$, $a\in\alg$, and $y\in\mcl{Y}$.
Here, for $x\in\mcl{Y}$, $U_{g_t(y,\cdot)}$ is the Koopman operator on $\lz$ with respect to the map $g_t(y,\cdot)$.
\end{definition}
\begin{remark}
The operator family $\{K_{\Phi_t}\}_{t\in\mathbb{R}}$ satisfies $K_{\Phi_s}K_{\Phi_t}=K_{\Phi_{s+t}}$ for any $s,t\in\mathbb{R}$ and $K_{\Phi_0}=I$.
However, it is not strongly continuous even for a simple case.
Let $\mcl{Z}=\mathbb{R}/2\pi\mathbb{Z}$ equipped with the normalized Haar measure on $\mcl{Z}$.
Let $g_t(y,z)=z+t\alpha$ for $\alpha\neq 0$.
For $v\equiv 1$ and $a=I$, we have
\begin{align*}
\Vert K_{\Phi_t}v\otimes a-v\otimes a\Vert_{\modu}^2
&=\bigg\Vert\int_{y\in\mcl{Y}}(U_{g_t(y,\cdot)}-I)^*(U_{g_t(y,\cdot)}-I) \mr{d}\mu(y)\bigg\Vert_{\alg}\\
&=\bigg\Vert\int_{y\in\mcl{Y}}(2I-U_{g_t(y,\cdot)}-U_{g_t(y,\cdot)}^*) \mr{d}\mu(y)\bigg\Vert_{\alg}\\
&=\Vert 2I-\tilde{U}^*M_{t}\tilde{U}-\tilde{U}^*M_{t}^*\tilde{U} \Vert_{\alg}\\
&=\Vert 2I-M_{t}-M_{t}^*\Vert_{\alg}\\
&=\sup_{n\in\mathbb{Z}}\vert 2-\mr{e}^{\sqrt{-1}nt\alpha}-\mr{e}^{-\sqrt{-1}nt\alpha}\vert,
\end{align*}
where $\tilde{U}:\lz\to L^2(\mathbb{Z})$ the unitary operator defined as $\gamma_i\mapsto e_i$, $\gamma_i(z)=\mr{e}^{\sqrt{-1}iz}$, and $e_i$ is the map on $\mathbb{Z}$ defined as $e_i(i)=1$ and $e_i(n)=0$ for $n\neq i$.
Moreover, $M_{t}:L^2(\mathbb{Z})\to L^2(\mathbb{Z})$ is the multiplication operator with respect to the map $n\mapsto \mr{e}^{\sqrt{-1}\alpha tn}$.
The third equality holds since
\begin{align*}
\tilde{U}^*M_{t}\tilde{U}\gamma_i=\tilde{U}^*\mr{e}^{\sqrt{-1}\alpha ti}e_i
=\mr{e}^{\sqrt{-1}\alpha ti}\gamma_i
=U_{g_t(y,\cdot)}\gamma_i.
\end{align*}
Let $\epsilon=\vert 2-\mr{e}^{\sqrt{-1}\alpha}-\mr{e}^{-\sqrt{-1}\alpha}\vert$.
For any $\delta>0$, let $n_0\in\mathbb{Z}$ such that $n_0\ge 1/\delta$ and let $t=1/n_0$.
Then, we have
\begin{align*}
\Vert K_{\Phi_t}v\otimes a-v\otimes a\Vert_{\modu}^2
&\ge \vert 2-\mr{e}^{\sqrt{-1}n_0t\alpha}-\mr{e}^{-\sqrt{-1}n_0t\alpha}\vert
=\epsilon.
\end{align*}
\end{remark}

We adopt the generator defined using a weaker topology than the topology of the Hilbert $C^*$-module.
\begin{definition}[Equicontinuous $C_0$-group~\cite{choe85}]\label{def:equicont_c0}
Let $M$ be a sequentially complete locally convex space and for any $t\in\mathbb{R}$, let $\kappa_t:M\to M$ be a linear operator on $M$ which satisfies
\begin{enumerate}
\item $\kappa_0=I$,
\item $\kappa_s\kappa_t=\kappa_{s+t}$ for any $s,t\in\mathbb{R}$,
 \item $\lim_{t\to 0}\kappa_tw=w$ for any $w\in M$,
 \item For any continuous seminorm $p$ on $M$, there exists a continuous seminorm $q$ such that $p(\kappa_tw)\le q(w)$ for any $w\in M$ and $t\in\mathbb{R}$.
\end{enumerate}
The family $\{\kappa_t\}_{t\in\mathbb{R}}$ is called an equicontinuous $C_0$-group.
\end{definition}
\begin{proposition}\label{prop:c_0group}
The space $\modu\subseteq \mcl{B}(\lz,\lx)$ equipped with the strong operator topology is a sequentially complete locally convex space.
In addition, assume $\mcl{Y}$ and $\mcl{Z}$ are locally compact Hausdorff spaces, $\mu$ and $\nu$ are regular probability measures, and $h$ and $g$ are continuous.
Then, $\{K_{\Phi_t}\}_{t\in\mathbb{R}}$ is an equicontinuous $C_0$-group.
\end{proposition}
To prove Proposition~\ref{prop:c_0group}, we use the following lemma:
\begin{lemma}\label{lem:continuous}
Let $\Omega$ and $\mcl{X}$ be topological spaces.
If a map $\Psi:\Omega\times\mcl{X}\to\mathbb{C}$ is continuous and compactly supported, then the map $\Omega\ni t\mapsto \Psi(t,\cdot)\in C_c(\mcl{X})$ is continuous.
Here, $C_c(\mcl{X})$ is the space of compactly supported continuous functions on $\mcl{X}$.
\end{lemma}
\begin{proof}
The statement follows from Lemma 4.16 by Eisner et al.~\cite{eisner16}.
\end{proof}
\begin{proof}[Proof of Proposition~\ref{prop:c_0group}]
\phantom{aaa}\\
({\bf $\modu$ is a sequentially complete locally convex space})
For $p\in\lz$, let $\Vert\cdot\Vert_p:\modu\to\mathbb{R}_+$ be defined as $\Vert w \Vert_p=\Vert wp\Vert_{\lx}$ for $w\in\modu\subseteq\mcl{B}(\lz,\lx)$.
Then, $\Vert\cdot\Vert_p$ is a seminorm in $\modu$.
Moreover, let $\{w_i\}_{i\in\mathbb{N}}$ be a countable Cauchy sequence in $\modu$.
Then, for any $v\in\lz$, $\{w_iv\}_{i\in\mathbb{N}}$ is a Cauchy sequence in the Hilbert space $\lx$.
Thus, there exists $\tilde{w}\in\lx$ such that $\lim_{i\to\infty}w_iv=\tilde{w}$.
Let $w:\lz\to \lx$ be the map defined as $w:v\mapsto \tilde{w}$.
Then, $w$ is linear and 
\begin{align*}
 \Vert wv\Vert_{\lx}=\Vert \lim_{i\to\infty}w_iv\Vert_{\lx}\le \sup_{i\in\mathbb{N}}\Vert w_iv\Vert_{\lx}
\le \sup_{i\in\mathbb{N}}\Vert w_i\Vert_{\modu}\,\Vert v\Vert_{\lz}
\end{align*}
for $v\in\lz$.
By the uniform boundedness principle, $\sup_{i\in\mathbb{N}}\Vert w_i\Vert<\infty$.
Thus, $w\in \mcl{B}(\lz,\lx)$.
Since $\modu\subseteq \mcl{B}(\lz,\lx)$ is closed with respect to the strong operator topology, we obtain $w\in\modu$.
Therefore, $\{w_i\}_{i\in\mathbb{N}}$ converges to $w$ in $\modu$.\vspace{.3cm}\\
({\bf $\{K_{\Phi_t}\}_{t\in\mathbb{R}}$ is an equicontinuous $C_0$-group})
For any $v\in\ly$, $a\in\alg$, and $u\in\lz$, we have
\begin{align*}
 \Vert (K_{\Phi_t}v\otimes a)u\Vert_{\lx}^2
&=\int_{y\in\mcl{Y}}\int_{z\in\mcl{Z}}\vert v(h_t(y))(au)(g_t(y,z))\vert^2\mr{d}\nu(z)\mr{d}\mu(y)\\
&=\int_{y\in\mcl{Y}}\int_{z\in\mcl{Z}}\vert v(y)(au)(z)\vert^2\mr{d}\nu(z)\mr{d}\mu(y)
=\Vert (v\otimes a)u\Vert_{\lx},
\end{align*}
which shows that the condition 4 of Definition~\ref{def:equicont_c0} is satisfied.

Regarding the condition 3, 
let $\epsilon>0$, let $\{\gamma_i\}_{i=1}^{\infty}$ be an orthonormal basis of $\lz$, and let $\mcl{D}=\{\sum_{i\in F}c_i\gamma_i\,\mid\,F\subset\mathbb{Z}:\mbox{ finite, }c_i\in\mathbb{C}\}$.
Since $C_c(\mcl{Z})$, $C_c(\mcl{Y})$, and $\mcl{D}$ are dense in $\lz$, $\ly$, and $\lz$, respectively, for any $i\in\mathbb{N}$ and any $v\in\ly$, $a\in\alg$, and $u\in\lz$, there exist $\tilde{v}\in C_c(\mcl{Y})$, $\tilde{\gamma}_i\in C_c(\mcl{Z})$, and $\tilde{u}\in\mcl{D}$ such that $\Vert \tilde{v}-v\Vert_{\ly}\le\epsilon$, $\Vert \tilde{\gamma}_i-a\gamma_i\Vert_{\lz}\le\epsilon/(\sqrt{2}^i)$, and $\Vert \tilde{u}-u\Vert_{\lz}\le\epsilon$.
Let $\tilde{a}=\sum_{i=1}^{\infty}\tilde{\gamma}_i\gamma_i'$, where the limit is taken with respect to the strong operator topology.
The operator $\tilde{a}$ is bounded since we have
\begin{align*}
&\Vert \tilde{a}u\Vert_{\lz}\le \Vert au\Vert_{\lz}+\Vert \tilde{a}u-au\Vert_{\lz}
\le\Vert a\Vert_{\alg} \Vert u\Vert_{\lz}+ \bigg\Vert \sum_{i=1}^{\infty}(a\gamma_i\gamma_i'u-\tilde{\gamma}_i\gamma_i'u)\bigg\Vert_{\lz}\\
&\qquad=\Vert a\Vert_{\alg} \Vert u\Vert_{\lz}+\bigg\Vert \sum_{i=1}^{\infty}(\gamma_i'u)(a\gamma_i-\tilde{\gamma}_i)\bigg\Vert_{\lz}\\
&\qquad\le\Vert a\Vert_{\alg} \Vert u\Vert_{\lz}+\bigg(\sum_{i=1}^{\infty}\vert\gamma_i'u\vert^2\bigg)^{1/2}\bigg(\sum_{i=1}^{\infty}\frac{\epsilon^2}{2^i}\bigg)^{1/2}
=\Vert u\Vert_{\lz}(\Vert a\Vert_{\alg}+\epsilon).
\end{align*}
Thus, we have $\tilde{a}\in\alg$.
In addition, we have
\begin{align*}
&\Vert (K_{\Phi_t}\tilde{v}\otimes \tilde{a})\tilde{u}-(\tilde{v}\otimes \tilde{a})\tilde{u}\Vert_{\lx}^2\\
&\qquad=\int_{y\in\mcl{Y}}\int_{z\in\mcl{Z}}\vert\tilde{v}(h_t(y))(U_{g_t(y,\cdot)}\tilde{a}\tilde{u})(z)-\tilde{v}(y)(\tilde{a}\tilde{u})(z)\vert^2\mr{d}\nu(z)\mr{d}\mu(y)\\
&\qquad=\int_{y\in\mcl{Y}}\int_{z\in\mcl{Z}}\vert\tilde{v}(h_t(y))(\tilde{a}\tilde{u})(g_t(y,z))-\tilde{v}(y)(\tilde{a}\tilde{u})(z)\vert^2\mr{d}\nu(z)\mr{d}\mu(y).
\end{align*}
Let $\Psi:\mathbb{R}\times \mcl{Y}\times\mcl{Z}\to\mathbb{C}$ be defined as $(t,y,z)\mapsto \tilde{v}(h_t(y))(\tilde{a}\tilde{u})(g_t(y,z))-\tilde{v}(y)(\tilde{a}\tilde{u})(z)$.
Since $\Psi$ is continuous, by Lemma~\ref{lem:continuous}, the map $t\mapsto\Psi(t,\cdot,\cdot)\in C_c(\mcl{X})$ is also continuous.
Thus, we have 
\begin{align*}
\lim_{t\to 0}\Vert (K_{\Phi_t}\tilde{v}\otimes \tilde{a})\tilde{u}-(\tilde{v}\otimes \tilde{a})\tilde{u}\Vert_{\lx}
&\le \lim_{t\to 0}\Vert (K_{\Phi_t}\tilde{v}\otimes \tilde{a})\tilde{u}-(\tilde{v}\otimes \tilde{a})\tilde{u}\Vert_{L^{\infty}(\mcl{X})}\\
&= \lim_{t\to 0}\Vert \Psi_t\Vert_{\infty}
=0,
\end{align*}
where $\Vert\cdot\Vert_{\infty}$ is the sup norm in $C_c(\mcl{X})$.
Therefore, $\lim_{t\to 0}\Vert (K_{\Phi_t}{v}\otimes {a}){u}-({v}\otimes {a}){u}\Vert_{\lx}=0$.
Indeed, we have
\begin{align*}
&\Vert (K_{\Phi_t}v\otimes a)u-(K_{\Phi_t}\tilde{v}\otimes \tilde{a})\tilde{u}\Vert_{\lx}
=\Vert (v\otimes a)u-(\tilde{v}\otimes \tilde{a})\tilde{u}\Vert_{\lx}\\
&\qquad\le\Vert (\tilde{v}\otimes (a-\tilde{a}))\tilde{u}\Vert_{\lx}
+\Vert ((v-\tilde{v})\otimes a)\tilde{u}\Vert_{\lx}+\Vert (v\otimes a)(u-\tilde{u})\Vert_{\lx}\\
&\qquad \le\Vert \tilde{v}\Vert_{\ly}\Vert(a-\tilde{a})\tilde{u} \Vert_{\lz}
+\Vert v-\tilde{v}\Vert_{\ly}\Vert {a}\tilde{u}\Vert_{\lz}
+\Vert v\Vert_{\ly}\Vert a\Vert_{\alg}\Vert u-\tilde{u}\Vert_{\lz}\\
&\qquad\le(\Vert v\Vert_{\ly}+\epsilon)\Vert u\Vert_{\lz}\epsilon+\epsilon\Vert au\Vert_{\lz}
+\Vert v\Vert_{\ly}\Vert a\Vert_{\alg}\epsilon.
\end{align*}
As a result, $\{K_{\Phi_t}\}_{t\in\mathbb{R}}$ satisfies the condition 3 of Definition~\ref{def:equicont_c0}.
\end{proof}
\begin{definition}
The generator $L_{\Phi}$ of $\{K_{\Phi_t}\}_{t\in\mathbb{R}}$ is defined as 
\begin{equation*}
L_{\Phi}w=\lim_{t\to 0}\frac{K_{\Phi_t}w-w}{t},
\end{equation*}
where the limit is with respect to the strong operator topology in $\modu$.
\end{definition}
\begin{proposition}[Choe, 1985~\cite{choe85}]
The generator $L_{\Phi}$ is a densely defined linear operator in $\modu$ with respect to the strong operator topology.
\end{proposition}

\subsection{Decomposition of $K_{\Phi_t}$ and $L_{\Phi}$}
We derive the eigenoperator decomposition for continuous systems.
In the following, we assume $\mcl{Y}$ and $\mcl{Z}$ are differentiable manifolds, $\mu$ and $\nu$ are regular probability measures, and $h$ and $g$ are differentiable.
\subsubsection{Fundamental decomposition}
We first define vectors to decompose the operator $K_{\Phi_t}$ using the cocycle.
\begin{definition}
For $s\in\mathbb{R}$, we define a linear operator $w_s:\lz\to\lx$ as
\begin{equation*}
(w_su)(y,z)=
(U_{g_s(y,\cdot)}u)(z).
\end{equation*}
\end{definition}

\begin{proposition}\label{prop:w_i_conti}
For $s\in\mathbb{R}$, we have $w_s\in\modu$.
Moreover, $K_{\Phi_t}w_s=w_{s+t}=w_s\cdot M_{s,t}$, where $M_{s,t}$ is a left $\alg$-linear multiplication operator on $\modu$ defined as $(w\cdot M_{s,t})(y)=w(y) U_{g_t(h_s(y),\cdot)}$.
\end{proposition}
\begin{proof}
We obtain $w_s\in\modu$ by Lemma~\ref{lem:winM}.
Moreover, we have
\begin{equation*}
K_{\Phi_t}w_s=U_{g_t(y,\cdot)}U_{g_s(h_t(y),\cdot)}=U_{g_s(h_t(y),g_t(y,\cdot))}=U_{g_{s+t}(y,\cdot)}=U_{g_s(y,\cdot)}U_{g_t(h_s(y),\cdot)}.
\end{equation*}
\end{proof}
\begin{proposition}\label{prop:decomposition_conti_fundamental}
For $s\in\mathbb{R}$ and $u\in C_c^1(\mcl{Z})$, let $\tilde{w}_{s,u}(y,z)=\frac{\partial u\circ g}{\partial t}(s,y,z)$.
Then, $(L_{\Phi}w_s)u=\tilde{w}_{s,u}$ and $L_{\Phi}w_s=w_s\cdot N_s$, where $(w\cdot N_s)(y)=w(y)(M_{\frac{\partial g}{\partial t}(0,h_s(y),\cdot)}\frac{\partial}{\partial z})$.
Here, $C_c^1(\mcl{Z})$ is the space of compactly supported and continuously differentiable functions on $\mcl{Z}$.
\end{proposition}
\begin{proof}
For $u\in C_c^1(\mcl{Z})$, we have
\begin{align*}
&\bigg\Vert \frac{1}{t}(K_{\Phi_t}w_s-w_s)u-\tilde{w}_{s,u}\bigg\Vert^2_{\lx}\\
&\qquad=\int_{y\in\mcl{Y}}\int_{z\in\mcl{Z}}\bigg\vert\frac{1}{t}(U_{g_t(y,\cdot)}U_{g_s(h_t(y),\cdot)}-U_{g_s(y,\cdot)})u(z)-\tilde{w}_{s,u}(y,z)\bigg\vert^2\mr{d}\nu(z)\mr{d}\mu(y)\\
&\qquad=\int_{y\in\mcl{Y}}\int_{z\in\mcl{Z}}\bigg\vert\frac{1}{t}(U_{g_{s+t}(y,\cdot)}-U_{g_s(y,\cdot)})u(z)-\tilde{w}_{s,u}(y,z)\bigg\vert^2\mr{d}\nu(z)\mr{d}\mu(y)\\
&\qquad=\int_{y\in\mcl{Y}}\int_{z\in\mcl{Z}}\bigg\vert\frac{1}{t}(u(g_{s+t}(y,z))-u(g_s(y,z)))-\tilde{w}_{s,u}(y,z)\bigg\vert^2\mr{d}\nu(z)\mr{d}\mu(y).
\end{align*}
Since $g$ is continuous, there exists $D>0$ such that for any $s\in\mathbb{R}$, $y\in\mcl{Y}$, and $z\in\mcl{Z}$, $\vert\frac{\partial u\circ g}{\partial t}(s,y,z)\vert <D$.
By the mean-value theorem, for any $y\in\mcl{Y}$ and $z\in\mcl{Z}$, there exists $c\in(s,s+t)$ for $t>0$ or $c\in (s+t,s)$ for $t<0$ such that 
\begin{equation*}
\bigg\vert\frac{1}{t}(u(g_{s+t}(y,z))-u(g_s(y,z)))\bigg\vert=\bigg\vert\frac{\partial u\circ g}{\partial t}(c,y,z)\bigg\vert\le D.
\end{equation*}
Thus, by the Lebesgue's dominated convergence theorem, we obtain
\begin{align*}
&\lim_{t\to 0}\bigg\Vert \frac{1}{t}(K_{\Phi_t}w_s-w_s)u-\tilde{w}_{s,u}\bigg\Vert^2_{\lx}\\
&\qquad=\int_{y\in\mcl{Y}}\int_{z\in\mcl{Z}}\lim_{t\to 0}\bigg\vert\frac{1}{t}(u(g_{s+t}(y,z))-u(g_s(y,z)))-\tilde{w}_{s,u}(y,z)\bigg\vert^2\mr{d}\nu(z)\mr{d}\mu(y)
=0.
\end{align*}
Thus, we have $(L_{\Phi}w_s)u=\tilde{w}_{s,u}$.
Moreover, $\tilde{w}_{s,u}$ is represented as 
\begin{align*}
\tilde{w}_{s,u}(y,z)&=\frac{\partial u\circ g}{\partial t}(s,y,z)
=\frac{\partial u}{\partial z}(g(s,y,z))\frac{\partial g}{\partial t}(s,y,z)\\
&=\bigg(U_{g_s(y,\cdot)}\frac{\partial u}{\partial z}\bigg)(z)\frac{\partial g}{\partial t}(s,y,z)
=U_{g_s(y,\cdot)}M_{\frac{\partial g}{\partial t}(s,y,g_s^{-1}(y,\cdot))}\frac{\partial }{\partial z}u(z).
\end{align*}
Since $g_s(y,g_{-s}(h_s(y),z))=g_s(h_{-s}(h_s(y)),g_{-s}(h_s(y),z))=g_0(h_s(y),z)=z$, $g_s(y,\cdot)^{-1}=g_{-s}(h_s(y),\cdot)$.
Thus, we have
\begin{align*}
\frac{\partial g}{\partial t}(s,y,g_s^{-1}(y,z))=\frac{\partial g}{\partial t}(s,h_{-s}(h_s(y)),g_{-s}(h_s(y),z))=\frac{\partial g}{\partial t}(0,h_s(y),z).
\end{align*}
\end{proof}

The vectors $w_s$ describe the dynamics on $\mcl{Z}$, which is specific for the skew product dynamical systems and we are interested in.

\begin{proposition}
The action of the Koopman operator $U_{\Phi_t}$ is decomposed into two parts as
\begin{equation*}
U_{\Phi_t}(v\otimes u)=U_{h_t}v\otimes U_{\Phi_s}u\circ g_{t-s}
\end{equation*}
for $v\in\ly$, $u\in\lz$, and $s,t\in\mathbb{R}$.
\end{proposition}
\begin{proof}
By the definition of $U_{\Phi_t}$, we have
\begin{align*}
U_{\Phi_t}(v\otimes u)(y,z)&=v(h_t(y))u(g_t(y,z))=v(h_t(y))u(g_{t-s}(h_s(y),g_s(y,z)))\\
&=U_{h_t}v(y)U_{\Phi_s}u\circ g_{t-s}(y,z).
\end{align*}
\end{proof}
Let 
\begin{equation*}
\mcl{W}_{0}=\bigg\{\sum_{s\in F}w_sc_s\,\mid\, F\subseteq\mathbb{R}:\ \mbox{finite set},\ c_s\in\alg\bigg\}
\end{equation*}
and $\mcl{W}$ be the completion of $\mcl{W}_{0}$ with respect to the norm in $\modu$ ($\mcl{W}$ is a submodule of $\modu$ and Hilbert $\alg$-module).
Moreover, for $s\in\mathbb{R}$ and $u\in\lz$, let $\tilde{w}_{u,s}\in\lx$ be defined as $\tilde{w}_{u,s}(y,z)=u(g_s(y,z))$.
Let
\begin{equation*}
\tilde{\mcl{W}}_{0}=\bigg\{\sum_{j=1}^n\sum_{s\in F}c_s\tilde{w}_{u_j,s}\,\mid\, n\in\mathbb{N},\ F\subseteq\mathbb{R}:\ \mbox{finite set},\ c_s\in\mathbb{C},\ u_j\in\lz\bigg\}
\end{equation*}
and $\tilde{\mcl{W}}$ be the completion of $\tilde{\mcl{W}}_{0}$ with respect to the norm in $\lx$.
\begin{proposition}
With the notation defined in Proposition~\ref{prop:UandK}, we have $K_{\Phi_t}\vert_{\mcl{W}}\,\iota_i\vert_{\tilde{\mcl{W}}}=\iota_i U_{\Phi_t}\vert_{\tilde{\mcl{W}}}$ for $t\in\mathbb{R}$ and $i=1,2,\ldots$.
\end{proposition}
\begin{proof}
For $u\in\lz$, $s\in\mathbb{R}$, and $i>0$, we have
\begin{align*}
(\iota_i \tilde{w}_{u,s})(y)=u(g_{s}(y,\cdot))\gamma_i'
=U_{g_s(y,\cdot)}u\gamma_i'
=w_s(y)(u\gamma_i').
\end{align*}
Thus, we obtain $\iota_i \tilde{w}_{u,s}\in\mcl{W}$.
Therefore, the range of $\iota_i\vert_{\tilde{\mcl{W}}}$ is contained in $\mcl{W}$.
The equality is deduced by the definitions of $K_{\Phi_t}$ and $U_{\Phi_t}$.
\end{proof}
\subsubsection{Further decomposition}
We further decompose $w_s$ and $N_{s}$ and obtain a more detailed decomposition of $L_{\Phi}\vert_{\mcl{W}}$.
For $y\in\mcl{Y}$, let $\mcl{V}_{1}(y),\mcl{V}_{2}(y),\ldots$ be a sequence of 
closed subspaces of $\lz$ which satisfies $\lz=\overline{\opn{Span}\{\bigcup_{j=1}^{\infty}\mcl{V}_{j}(y)\}}$ for a.s. $y\in\mcl{Y}$.
For $s\in\mathbb{R}$ and $j=1,2,\ldots$, we define a linear map $\hat{w}_{s,j}$ from $\lz$ to $\lx$ as $(\hat{w}_{s,j}u)(y,z)=(w_s(y)p_{j}(h_s(y)) u)(y,z)$, where $p_{j}(y):\lz\to\mcl{V}_{j}(y)$ is the projection onto $\mcl{V}_{j}(y)$.
Assume $p_{j}(y)$ satisfies $(y,z)\mapsto (p_{j}(y)u)(z)\in\lx$.
We denote by $p_{j}$ the linear operator from $\lz\to\lx$ defined as $p_{j}u(y,z)=(p_{j}(y)u)(z)$.
For each $j=1,2,\ldots$, the following proposition holds.
Here, we define a differential operator $V_{\Phi}$ by
\begin{align}
V_{\Phi}v(y,z)=\frac{\partial v}{\partial y}(y,z)\frac{\partial h}{\partial t}(0,y)+\frac{\partial v}{\partial z}(y,z)\frac{\partial g}{\partial t}(0,y,z)\label{eq:generator}
\end{align}
for $v\in C^1_c(\mcl{X})$.

\begin{theorem}[Eigenoperator decomposition for continuous-time systems]\label{prop:hat_w_conti}
For $s\in\mathbb{R}$ and $u\in\lz$, let 
\begin{align*}
\tilde{w}_{s,u,j}(y,z)=\frac{\partial (p_{j}(h_t(y))u\circ g(t,y,z))}{\partial t}\bigg\vert_{t=s}.
\end{align*}
Assume for any $u\in p_{j}^{-1}(C_c^1(\mcl{X}))$ and any $y\in\mcl{Y}$, $\frac{\partial(p_{j}(h_t(y))u)(g_t(y,\cdot))}{\partial t}\big\vert_{t=0}=(V_{\Phi}p_{j}u)(y,\cdot)\in\mcl{V}_{j}(y)$. 
Then, $(L_{\Phi}\hat{w}_{s,j})u=\tilde{w}_{s,u,j}=(\hat{w}_{s,j}\cdot \hat{N}_{s,j})u$, where $\hat{N}_{s,j}$ is defined as $(w\cdot \hat{N}_{s,j})u(y,z)=w(y)(V_{\Phi}p_{j}u)(h_s(y),z)$.
\end{theorem}
We call $\hat{N}_{s,j}$ an eigenoperator and $\hat{w}_{s,j}$ an eigenvector.

\begin{proof}
For $u\in p_{j}^{-1}(C_c^1(\mcl{X}))$, we have
\begin{align*}
&\bigg\Vert \frac{1}{t}(K_{\Phi_t}\hat{w}_{s,j}-\hat{w}_{s,j})u-\tilde{w}_{s,u,j}\bigg\Vert^2_{\lx}\\
&=\int_{y\in\mcl{Y}}\int_{z\in\mcl{Z}}\bigg\vert\frac{1}{t}(U_{g_t(y,\cdot)}U_{g_s(h_t(y),\cdot)}p_{j}(h_{s+t}(y))-U_{g_s(y,\cdot)}p_{j}(h_s(y)))u(z)-\tilde{w}_{s,u,j}(y,z)\bigg\vert^2\mr{d}\nu(z)\mr{d}\mu(y)\\
&=\int_{y\in\mcl{Y}}\int_{z\in\mcl{Z}}\bigg\vert\frac{1}{t}(U_{g_{s+t}(y,\cdot)}p_{j}(h_{s+t}(y))-U_{g_s(y,\cdot)}p_{j}(h_s(y)))u(z)-\tilde{w}_{s,u,j}(y,z)\bigg\vert^2\mr{d}\nu(z)\mr{d}\mu(y)\\
&=\int_{y\in\mcl{Y}}\int_{z\in\mcl{Z}}\bigg\vert\frac{1}{t}((p_{j}(h_{s+t}(y))u)(g_{s+t}(y,z))-(p_j(h_s(y))u)(g_s(y,z)))-\tilde{w}_{s,u,j}(y,z)\bigg\vert^2\mr{d}\nu(z)\mr{d}\mu(y)\\
&\to 0\ (t\to 0).
\end{align*}
Thus, we have $(L_{\Phi}\hat{w}_{s,j})u=\tilde{w}_{s,u,j}$.
Moreover, $\tilde{w}_{s,u,j}$ is represented as 
\begin{align*}
\tilde{w}_{s,u,j}(y,z)&=\frac{\partial}{\partial t} (U_{g_t(y,\cdot)}p_{j}(h_t(y))u)(z)\bigg\vert_{t=s}
=U_{g_s(y,\cdot)}\frac{\partial}{\partial t} (U_{g_t(h_s(y),\cdot)}p_{j}(h_{s+t}(y))u)(z)\bigg\vert_{t=0}\\
&=U_{g_s(y,\cdot)}\frac{\partial (p_{j}(h_{t}(h_s(y)))u)(g(t,h_s(y),z))}{\partial t}\bigg\vert_{t=0}\\
&=U_{g_s(y,\cdot)}p_{j}(h_s(y))\frac{\partial (p_{s,j}(h_t(y))u)(g(t,h_s(y),z))}{\partial t}\bigg\vert_{t=0}\\
&=\hat{w}_{s,j}(y)\bigg(\frac{\partial (p_{s,j}(y)u)(z)}{\partial y}\bigg\vert_{\substack{y=h(0,y),\\z=g(0,h_s(y),z)}}\frac{\partial h(t,y)}{\partial t}\bigg\vert_{t=0}\\
&\qquad\qquad\qquad+\frac{\partial (p_{s,j}(y)u)(z)}{\partial z}\bigg\vert_{\substack{y=h(0,y),\\z=g(0,h_s(y),z)}}\frac{\partial g(t,h_s(y),z)}{\partial t}\bigg\vert_{t=0}\bigg)\\
&=\hat{w}_{s,j}(y)\bigg(\frac{\partial}{\partial y}\cdot\frac{\partial h}{\partial t}(0,y)+\frac{\partial }{\partial z}\cdot \frac{\partial g}{\partial t}(0,h_s(y),z)\bigg) p_{s,j}u(y,z),
\end{align*}
where $p_{s,j}(y)=p_j(h_s(y))$.
Furthermore, we have
\begin{align*}
&\bigg(\frac{\partial}{\partial y}\cdot \frac{\partial h}{\partial t}(0,y)+\frac{\partial }{\partial z}\cdot\frac{\partial g}{\partial t}(0,h_s(y),z)\bigg) p_{s,j}u(y,z)\\
&\qquad=\frac{\partial p_{s,j}u}{\partial y}(y,z)\frac{\partial h}{\partial t}(0,y)+\frac{\partial p_{s,j}u}{\partial z}(y,z)\frac{\partial g}{\partial t}(0,h_s(y),z)\\
&\qquad=\frac{\partial p_{j}u}{\partial y}(h_s(y),z)\frac{\partial h_s}{\partial y}(y)\frac{\partial h}{\partial t}(0,y)+\frac{\partial p_{j}u}{\partial z}(h_s(y),z)\frac{\partial g}{\partial t}(0,h_s(y),z)\\
&\qquad=\frac{\partial p_{j}u}{\partial y}(h_s(y),z)\frac{\partial h_s(h(t,y))}{\partial t}\bigg\vert_{t=0}
+\frac{\partial p_{j}u}{\partial z}(h_s(y),z)\frac{\partial g}{\partial t}(0,h_s(y),z)\\
&\qquad=\frac{\partial p_{j}u}{\partial y}(h_s(y),z)\frac{\partial h}{\partial t}(0,h_s(y))
+\frac{\partial p_{j}u}{\partial z}(h_s(y),z)\frac{\partial g}{\partial t}(0,h_s(y),z)
=V_{\Phi}p_ju(h_s(y),z).
\end{align*}
\end{proof}

By Proposition~\ref{prop:spectrum_finitedim}, we have the following proposition regarding the spectrum of $\hat{N}_{s,j}$.
\begin{proposition}\label{prop:eigenop_finite}
Assume $\opn{dim}(\mcl{V}_j(y))$ is finite and constant with respect to $y \in \mcl{Y}$.
Then, we have $\sigma(\hat{N}_{0,j})=\sigma(\hat{N}_{s,j})$ for any $s\in\mathbb{R}$.
\end{proposition}
\begin{proof}
Since $h_s$ is measure preserving for any $s\in\mathbb{R}$, by Proposition~\ref{prop:spectrum_finitedim}, we have
\begin{align*}
\sigma(\hat{N}_{0,j})
&=\{\lambda\in\mathbb{C}\,\mid\,^{\forall}\epsilon>0,\ \mu(\{y\in\mcl{Y}\,\mid\,\lambda\in\sigma_{\epsilon}(V_{\Phi}p_j(y))\})>0\}\\
&=\{\lambda\in\mathbb{C}\,\mid\,^{\forall}\epsilon>0,\ \mu(\{y\in\mcl{Y}\,\mid\,\lambda\in\sigma_{\epsilon}(V_{\Phi}p_j(h_s(y)))\})>0\}
=\sigma(\hat{N}_{s,j}).
\end{align*}
\end{proof}

\begin{remark}\label{rem:assumption}
If $U_{g_t(y,\cdot)}\mcl{V}_{j}(h_t(y))\subseteq\mcl{V}(y)$ for any $t$ in a neighborhood of $0$ in $\mathbb{R}$, then we have  
\begin{align*}
&\frac{\partial}{\partial t} (U_{g_t(y,\cdot)}p_{j}(h_t(y))u)(z)\bigg\vert_{t=0}
=\lim_{t\to 0}\frac{1}{t}(U_{g_t(y,\cdot)}p_{j}(h_t(y))-p_{j}(y))u(z)\\
&\quad=\lim_{t\to 0}\frac{1}{t}p_{j}(y)(U_{g_t(y,\cdot)}p_{j}(h_t(y))-p_{j}(y))u(z)
=p_{j}(y)\frac{\partial}{\partial t} (U_{g_t(y,\cdot)}p_{j}(h_t(y))u)(z)\bigg\vert_{t=0}.
\end{align*}
Thus, the assumption $\frac{\partial(p_{j}(h_t(y))u)(g_t(y,\cdot))}{\partial t}\big\vert_{t=0}\in\mcl{V}_{j}(y)$ in Theorem~\ref{prop:hat_w_conti} is satisfied.
\end{remark}

\begin{example}\label{ex:rotation_conti_spectrum}
Let $\mcl{Y}=\mcl{Z}=\mathbb{R}/2\pi\mathbb{Z}=:\mathbb{T}$.
For $\alpha,\beta>0$, consider the following continuous dynamical system:
\begin{equation}
\bigg(\frac{\mr{d}y(t)}{\mr{d}t},\frac{\mr{d}z(t)}{\mr{d}t}\bigg)=(1,\alpha(1+\beta\cos(y(t)))).\label{eq:rotation_conti}
\end{equation}
In this case, we have $\frac{\partial g}{\partial t}(0,y,z)=\alpha(1+\beta\cos(y(0)))=\alpha(1+\beta\cos(y))$ and $h_s(y)=y+s$.
Let $\gamma_{k,j}(y,z)=\mr{e}^{\sqrt{-1}(ky+jz)}$.
Then, we have
\begin{align*}
V_{\Phi}\gamma_{k,j}(y,z)=\big(\sqrt{-1}k+\sqrt{-1}j\alpha(1+\beta\cos(y))\big)\gamma_{k,j}(y,z).
\end{align*}
Let $\mcl{V}_{j}=\overline{\opn{Span}\{\gamma_{k,j}\,\mid\,k\in\mathbb{Z}\}}$.
We can see $\mcl{V}_{j}$ is an invariant subspace of $V_{\Phi}$.
In addition, let $\mcl{V}_{j}(y)=\overline{\opn{Span}\{\gamma_{k,j}(y,\cdot)\,\mid\,k\in\mathbb{Z}\}}$, and let $p_{j}(y)$ be the projection onto $\mcl{V}_{j}(y)$.
Then, since we have
\begin{align*}
(V_{\Phi}p_j)(y)\gamma_{k,j}(y,\cdot)=(V_{\Phi}\gamma_{k,j})(y,\cdot),
\end{align*}
the spectrum of $(V_{\Phi}p_{j})(y)$ is calculated as
\begin{align*}
\sigma((V_{\Phi}p_{j})(y))=\{\sqrt{-1}k+\sqrt{-1}j\alpha(1+\beta\cos(y))\,\mid\,k\in\mathbb{Z}\}.
\end{align*}
Therefore, we have 
\begin{align*}
\bigcup_{y\in\mcl{Y}}\sigma((V_{\Phi}p_{j})(y))=
\bigcup_{y\in\mcl{Y}}\sigma((V_{\Phi}p_{j})(h_s(y)))=
\bigcup_{y\in\mcl{Y}}\{\sqrt{-1}k+\sqrt{-1}j\alpha(1+\beta\cos(y))\,\mid\,k\in\mathbb{Z}\}.
\end{align*}
Regarding $\hat{w}_{s,j}$, we have
\begin{align*}
U_{g_s(y,\cdot)}\gamma_{k,j}(y,\cdot)&=\gamma_{k,j}(y,\cdot+\alpha(s+\beta(\sin(y+s)-\sin(y))))\\
&=\gamma_{k,j}(y,\cdot)\mr{e}^{\sqrt{-1}j\alpha(s+\beta(\sin(y+s)-\sin(y)))}.
\end{align*}
Thus, the spectrum of the family of operators $\{\hat{w}_{s,j}(y)\}_s$ on $\lz$ is $\gamma_j\big(\alpha s+\alpha\beta(\sin(y+s)-\sin(y))\big)$.
\end{example}
In the following subsections, we will generalize the arguments in Example~\ref{ex:rotation_conti_spectrum}. 

\subsubsection{Construction of the generalized Oseledets space $\mcl{V}_{j}(y)$ using a function space on $\mcl{X}$}
We show how we can construct the generalized Oseledets space $\mcl{V}_{j}(y)$ required for obtaining $p_j$ appearing in Theorem~\ref{prop:hat_w_conti}.
In this subsection, we assume $\mcl{Y}$ is compact.
Let 
\begin{align}
\moduu=C(\mcl{Y})\otimes\lz\label{eq:N_def}
\end{align}
be the Hilbert $C(\mcl{Y})$-module.
Note that a Hilbert $C^*$-module is also a Banach space.
Here, we just regard $\mcl{N}$ as a Banach space equipped with the norm $\Vert a\otimes u\Vert_{\moduu}^2= \sup_{y\in\mcl{Y}}\int_{z\in\mcl{Z}}\vert a(y)u(z)\vert^2 \mr{d}\nu(z)$. 
For $t\in\mathbb{R}$, let $U_{\Phi_{t}}$ be the Koopman operator on $\mcl{N}$ with respect to $\Phi_{t}$.
\begin{proposition}
The family of operators $\{U_{\Phi_{t}}\}_{t\in\mathbb{R}}$ is a strongly continuous one-parameter group.
\end{proposition}
\begin{proof}
Let $v\in C(\mcl{Y})\otimes_{\opn{alg}} C_c(\mcl{Z})$ and let $\epsilon>0$.
Then, there exists $\delta>0$ such that for any $\vert t\vert\le\delta$, $y\in\mcl{Y}$, and $z\in\mcl{Z}$, $\vert v(h_t(y),g_t(y,z))-v(y,z)\vert\le \epsilon$.
Thus, we have
\begin{equation}
\Vert U_{\Phi_{t}}v-v\Vert_{\moduu}
=\sup_{y\in\mcl{Y}}\bigg( \int_{z\in\mcl{Z}}\vert v(h_t(y),g_t(y,z))-v(y,z)\vert^2\mr{d}\nu(z)\bigg)^{1/2}
\le \epsilon.\label{eq:stongly_sontinuous}
\end{equation}
In addition, for any $v\in C(\mcl{Y})\otimes_{\opn{alg}} \lz$, we have
\begin{align*}
\Vert U_{\Phi_{t}}v\Vert_{\moduu}
&=\sup_{y\in\mcl{Y}}\bigg( \int_{z\in\mcl{Z}}\vert v(h_t(y),g_t(y,z))\vert^2\mr{d}\nu(z)\bigg)^{1/2}\\
&=\sup_{y\in\mcl{Y}}\bigg( \int_{z\in\mcl{Z}}\vert v(y,z)\vert^2\mr{d}\nu(z)\bigg)^{1/2}
=\Vert v\Vert_{\moduu}.
\end{align*}
Since $C(\mcl{Y})\otimes_{\opn{alg}} C_c(\mcl{Z})$ is dense in $\moduu$, Eq.~\eqref{eq:stongly_sontinuous} is satisfied for any $v\in \moduu$.
\end{proof}
We note that the generator of $\{U_{\Phi_{t}}\}_{t\in\mathbb{R}}$ is $V_{\Phi}$ defined in Eq.~\eqref{eq:generator}.
%
If we set $\mcl{V}_j$ as $\mcl{V}$ in the following proposition, it satisfies the assumption of Theorem~\ref{prop:hat_w_conti} (see also Remark~\ref{rem:assumption}).

\begin{proposition}\label{prop:construction_invariant_sp}
Let $\mcl{V}$ be an invariant subspace of $U_{\Phi_t}$ and let $\mcl{V}(y)=\overline{R_y\mcl{V}}$.
Then, we have $U_{g_t(y,\cdot)}\mcl{V}(h_t(y))\subseteq \mcl{V}(y)$.
Here, $R_{y}:\moduu\to\lz$ be a linear map defined as $R_{y}(a\otimes u)(z)=a(y)u(z)$ for $y\in\mcl{Y}$.
\end{proposition}
\begin{proof}
For $t\in\mathbb{R}$, $y\in\mcl{Y}$, and $v\in C(\mcl{Y})\otimes_{\opn{alg}}\lz$, we have
\begin{align*}
U_{g_t(y,\cdot)}R_{h_t(y)}v(z)
=v(h_t(y),g_t(y,z))
=R_{y}U_{\Phi_t}v(z).
\end{align*}
Thus, we have
$U_{g_t(y,\cdot)}R_{h_t(y)}=R_{y}U_{\Phi_{t}}$.
Since $\mcl{V}$ is an invariant subspace of $U_{\Phi_t}$, we have $U_{g_t(y,\cdot)}\mcl{V}(h_t(y))\subseteq \mcl{V}(y)$.
\end{proof}

The following proposition shows an example of $\mcl{V}$ constructed in Proposition~\ref{prop:construction_invariant_sp}.
It is for a simple case where $V_{\Phi}$ has an eigenvalue, but provides us with an intuition of what the eigenoperators describe.
\begin{proposition}\label{prop:simple_case}
Assume there exists $y\in\mcl{Y}$ such that $\{h(t,y)\,\mid\,t\in\mathbb{R}\}$ is dense in $\mcl{Y}$.
Assume $V_{\Phi}$ has an eigenvalue $\lambda$ and the corresponding eigenvector $\tilde{v}\in C_c^1(\mcl{X})$.
Then, there exists $C\ge 0$ such that for a.s. $y\in\mcl{Y}$, $\Vert \tilde{v}(y,\cdot)\Vert_{\lz}=C$.
Assume $C>0$ and let $p(y)u={v(y,\cdot)v(y,\cdot)^*u}$ for $u\in\lz$, where $v=\tilde{v}/C$.
Then, $p(y)(V_{\Phi}pu)(y,\cdot)=(V_{\Phi}pu)(y,\cdot)$ and 
\begin{align*}
\sigma((V_{\Phi}p)(y))=\lambda-\int_{\mcl{Z}}\frac{\partial {v}}{\partial y}(y,z)\overline{v}(y,z)\mr{d}\nu(z)\frac{\partial h}{\partial t}(0,y)=\int_{\mcl{Z}}\frac{\partial {v}}{\partial z}(y,z)\frac{\partial {g}}{\partial t}(0,y,z)\overline{v}(y,z)\mr{d}\nu(z).
\end{align*}
Moreover, 
$\sigma((V_{\Phi}p)(y))\subseteq \sqrt{-1}\mathbb{R}$.
\end{proposition}
\begin{proof}
The vector $\tilde{v}$ is an eigenvector of $U_{\Phi_t}$ for any $t\in\mathbb{R}$, and its corresponding eigenvalue is $\mr{e}^{\lambda t}$ ($\lambda\in\sqrt{-1}\mathbb{R}$).
Thus, we have
\begin{align*}
\int_{\mcl{Z}}\tilde{v}(y,z)\overline{\tilde{v}}(y,z)\mr{d}\nu(z)
&=\int_{\mcl{Z}}\mr{e}^{-\lambda t}U_{\Phi_t}\tilde{v}(y,z)\overline{\mr{e}^{-\lambda t}U_{\Phi_t}\tilde{v}(y,z)}\mr{d}\nu(z)\\
&=\int_{\mcl{Z}}\tilde{v}(h_t(y),g_t(y,z))\overline{\tilde{v}}(h_t(y),g_t(y,z))\mr{d}\nu(z)\\
&=\int_{\mcl{Z}}\tilde{v}(h_t(y),z)\overline{\tilde{v}}(h_t(y),z)\mr{d}\nu(z).
\end{align*}

For $u\in\lz$, we have
\begin{align*}
&(V_{\Phi}pu)(y,z)\\
&=\bigg(\frac{\partial v}{\partial y}(y,z)\int_{\mcl{Z}}\overline{v}(y,z)u(z)\mr{d}\nu(z)+v(y,z)\int_{\mcl{Z}}\frac{\partial \overline{v}}{\partial y}(y,z)u(z)\mr{d}\nu(z)\bigg)\frac{\partial h}{\partial t}(0,y)\\
&\qquad\qquad+\frac{\partial v}{\partial z}(y,z)\int_{\mcl{Z}}\overline{v}(y,z)u(z)\mr{d}\nu(z)\frac{\partial g}{\partial t}(0,y,z)\\
&=(V_{\Phi}v)(y,z)\int_{\mcl{Z}}\overline{v}(y,z)u(z)\mr{d}\nu(z)+v(y,z)\int_{\mcl{Z}}\frac{\partial \overline{v}}{\partial y}(y,z)u(z)\mr{d}\nu(z)\frac{\partial h}{\partial t}(0,y)\\
&=\lambda v(y,z)\int_{\mcl{Z}}\overline{v}(y,z)u(z)\mr{d}\nu(z)+v(y,z)\int_{\mcl{Z}}\frac{\partial \overline{v}}{\partial y}(y,z)u(z)\mr{d}\nu(z)\frac{\partial h}{\partial t}(0,y).
\end{align*}
Thus, $p(y)(V_{\Phi}pu)(y)=(V_{\Phi}pu)(y)$.
In addition, we have
\begin{align*}
(V_{\Phi}p)(y)v(y,\cdot)&=\lambda v(y,\cdot)\int_{\mcl{Z}}\overline{v}(y,z)v(y,z)\mr{d}\nu(z)+
v(y,\cdot)\int_{\mcl{Z}}\frac{\partial \overline{v}}{\partial y}(y,z)v(y,z)\mr{d}\nu(z)\frac{\partial h}{\partial t}(0,y)\\
&=\bigg(\lambda+\int_{\mcl{Z}}\frac{\partial \overline{v}}{\partial y}(y,z)v(y,z)\mr{d}\nu(z)\frac{\partial h}{\partial t}(0,y)\bigg)v(y,\cdot).
\end{align*}
Moreover, we have
\begin{align*}
0&=\int_{\mcl{Z}}\frac{\partial v(y,z)\overline{v}(y,z)}{\partial y}\mr{d}\nu(z)
=\int_{\mcl{Z}}\frac{\partial \overline{v}}{\partial y}(y,z)v(y,z)\mr{d}\nu(z)
+\int_{\mcl{Z}}\frac{\partial {v}}{\partial y}(y,z)\overline{v}(y,z)\mr{d}\nu(z)\\
&=2\Re\bigg(\int_{\mcl{Z}}\frac{\partial \overline{v}}{\partial y}(y,z)v(y,z)\mr{d}\nu(z)\bigg),
\end{align*}
and
\begin{align*}
\int_{\mcl{Z}}\frac{\partial \overline{v}}{\partial y}(y,z)v(y,z)\mr{d}\nu(z)\frac{\partial h}{\partial t}(0,y)
&=-\int_{\mcl{Z}}\frac{\partial {v}}{\partial y}(y,z)\overline{v}(y,z)\mr{d}\nu(z)\frac{\partial h}{\partial t}(0,y)\\
&=-\int_{\mcl{Z}}\bigg(V_{\Phi}v(y,z)-\frac{\partial {v}}{\partial z}(y,z)\frac{\partial g}{\partial t}(0,y,z)\bigg)\overline{v}(y,z)\mr{d}\nu(z)\\
&=\int_{\mcl{Z}}\bigg(-\lambda v(y,z)+\frac{\partial {v}}{\partial z}(y,z)\frac{\partial g}{\partial t}(0,y,z)\bigg)\overline{v}(y,z)\mr{d}\nu(z).
\end{align*}
\end{proof}

\begin{remark}
Proposition~\ref{prop:simple_case} implies that the eigenoperators have the information of the dynamics on $\mcl{Z}$ for each $y$, which cannot be extracted by eigenvalues of $V_{\Phi}$.
Indeed, if an eigenvector $v_j$ of $V_{\Phi}$ depends only on $y$, then $\sigma(\hat{N}_{s,j})=\sigma(V_{\Phi}p(y))=0$.
On the other hand, the corresponding eigenvalue of $V_{\Phi}$ can be nonzero.
\end{remark}

\subsubsection{Approximation of $V_{\Phi}$ in RKHS}
For numerical computations, to construct the subspace $\mcl{V}_j(y)$, we need to approximate $V_{\Phi}$ using RKHSs.
Approximating the generator of the Koopman operator in RKHSs was proposed by Das et al.~\cite{das21}.
Here, we apply a similar technique to approximating $V_{\Phi}$.
In this subsection, we assume $\mcl{Y}=\mathbb{T}$.
We also assume $\mcl{Z}$ is compact and $\nu$ is a Borel probability measure satisfying $\opn{supp}(\nu)=\mcl{Z}$.

Let $\phi_i(l)=\mr{e}^{\sqrt{-1}il}$ and $\lambda_i=\mr{e}^{-\vert i\vert}$ for $i\in\mathbb{Z}$ and $l\in\mathbb{T}$.
Let $p_1:\mathbb{T}\times\mathbb{T}\to\mathbb{C}$ be the positive definite kernel defined as $p_1(l_1,l_2)=\sum_{i\in\mathbb{Z}}\lambda_i\phi_i(l_1)\overline{\phi_i(l_2)}$.
In addition, let $p_2:\mcl{Z}\times\mcl{Z}\to \mathbb{C}$ be a positive definite kernel, and let $\tilde{P}:\lz\to\lz$ be the integral operator with respect to $p_2$.
Let $\tilde{\lambda}_1\ge\tilde{\lambda}_2\ge\ldots> 0$ and $\tilde{\phi}_1,\tilde{\phi}_2,\ldots$ be eigenvalues and the corresponding orthonormal eigenvectors of $\tilde{P}$, respectively.
By Mercer's theorem, $p_2(z_1,z_2)=\sum_{i=1}^{\infty}\tilde{\lambda}_i\tilde{\phi}_i(z_1)\overline{\tilde{\phi}_i(z_2)}$, where the sum converges uniformly on $\mcl{Z}\times\mcl{Z}$.
Let $\tau>0$ and let $\lambda_{\tau,i,j}=\mr{e}^{\tau(1-\lambda_i^{-1}\tilde{\lambda}_j^{-1})}$ for $i\in\mathbb{Z}$ and $j=1,2,\ldots$.
Let 
\begin{align*}
p_{\tau}((l_1,z_1),(l_2,z_2))=\sum_{i\in\mathbb{Z}}\sum_{j=1}^{\infty}\lambda_{\tau,i,j}\phi_i(l_1)\tilde{\phi}_j(z_1)\overline{\tilde{\phi}_j(z_2)}\overline{\phi_i(l_2)}
\end{align*}
and let $\hil_{\tau}$ be the RKHS associated with $p_{\tau}$.
In addition, let ${P}_{\tau}$ be the integral operator with respect to $p_{\tau}$.

\begin{proposition}\label{lem:Ptau}
Let $\iota_{\tau}:\hil_{\tau}\to\moduu$ be the inclusion map, where $\moduu$ is defined as Eq.~\eqref{eq:N_def}.
Then, for any $v\in \mcl{N}$, $\Vert \iota_{\tau}P_{\tau}v-v\Vert_{\moduu}$ converges to $0$ as $\tau\to 0$.
\end{proposition}
\begin{proof}
Let $\psi_{i,j}=\phi_i\otimes\tilde{\phi}_j$.
Since $\{\tilde{\phi}_i\}_{i=1}^{\infty}$ is an orthonormal basis in $\lz$, the subspace $\{\sum_{i=-n}^n\sum_{j=1}^ma_{i,j}\psi_{i,j}\,\mid\,n,m\in\mathbb{N},a_{i,j}\in\mathbb{C}\}$ is dense in $\moduu$ with $\Vert\psi_{i,j}\Vert_{\moduu}=1$.
In addition, since we have $\iota_{\tau}P_{\tau}\psi_{i,j}=\lambda_{\tau,i,j}\psi_{i,j}$ and $0\le \lambda_{\tau,i,j}\le 1$, we obtain $\Vert \iota_{\tau}P_{\tau}\Vert_{\moduu}\le 1$.
For any $\epsilon>0$ and $v\in\moduu$, there exist $n,m\in\mathbb{N}$, $a_{i,j}\in \mathbb{C}$, and $\tau_0>0$ such that $\Vert \sum_{i=-n}^n\sum_{j=1}^ma_{i,j}\psi_{i,j}-v\Vert_{\moduu}\le \epsilon$ and $(1-\lambda_{\tau,i,j})(\sum_{i=-n}^n\sum_{j=1}^m\vert a_{i,j}\vert)\le \epsilon$ for $i=-n,\ldots,n$, $j=1,\ldots,m$, and $\tau\le\tau_0$.
Thus, for $\tau\le\tau_0$, we have
\begin{align*}
&\Vert \iota_{\tau}P_{\tau}v-v\Vert_{\moduu}\\
&\le\bigg\Vert \iota_{\tau}P_{\tau}v- \iota_{\tau}P_{\tau}\sum_{i=-n}^n\sum_{j=1}^ma_{i,j}\psi_{i,j} \bigg\Vert_{\moduu}
+\bigg\Vert \iota_{\tau}P_{\tau}\sum_{i=-n}^n\sum_{j=1}^ma_{i,j}\psi_{i,j} -\sum_{i=-n}^n\sum_{j=1}^ma_{i,j}\psi_{i,j}\bigg\Vert_{\moduu}\\
&\qquad\qquad+\bigg\Vert \sum_{i=-n}^n\sum_{j=1}^ma_{i,j}\psi_{i,j}-v\bigg\Vert_{\moduu}\\
&\le\epsilon+\sum_{i=-n}^n\sum_{j=1}^m(1-\lambda_{\tau,i,j})\vert a_{i,j}\vert+\epsilon=3\epsilon.
\end{align*}
\end{proof}

By Proposition~\ref{lem:Ptau}, we can see that $V_{\Phi}$ can be approximated by $P_{\tau}V_{\Phi}\iota_{\tau}$ in the following sense.
\begin{corollary}
Let $\tau_0>0$.
For any $v\in\hil_{\tau_0}$, $\Vert \iota_{\tau}P_{\tau}V_{\Phi}\iota_{\tau}v-V_{\Phi}\iota_{\tau}v\Vert_{\moduu}$ converges to $0$ as $\tau\to 0$.
\end{corollary}

\section{Numerical examples}
\label{sec:continuous_time_examples}
We numerically investigate the eigenoperator decomposition.
\subsection{Moving Gaussian vortex}
\begin{figure}
    \begin{tabularx}{\textwidth}{YYYYY}
    \includegraphics[scale=0.23]{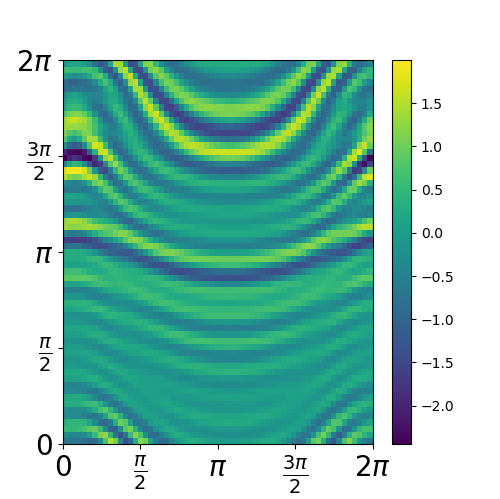}\newline $d=1$&
    \includegraphics[scale=0.23]{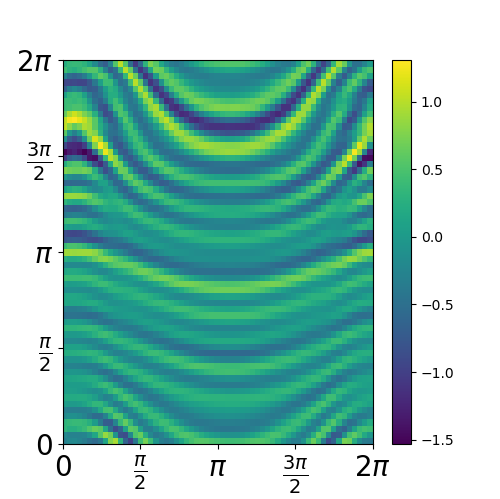}\newline $d=2$&
    \includegraphics[scale=0.23]{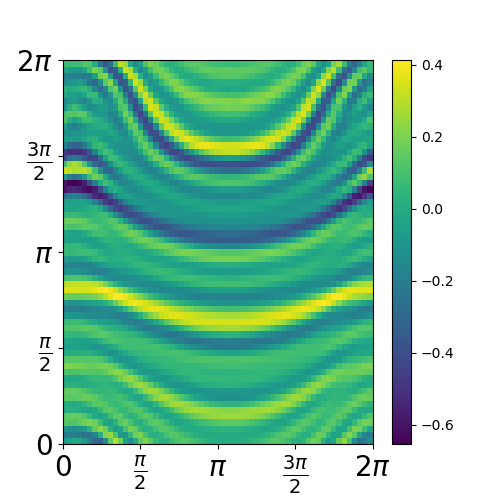}\newline $d=10$&
    \includegraphics[scale=0.23]{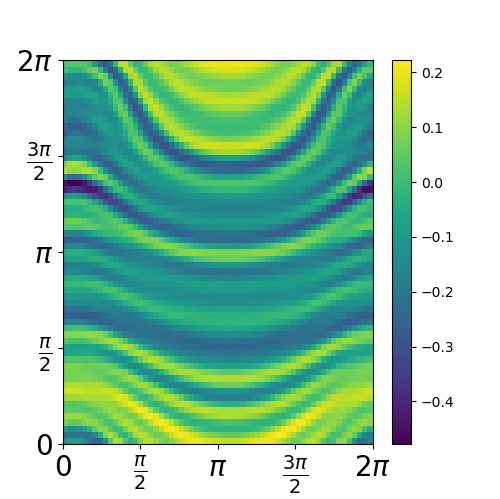}\newline $d=20$&
    \includegraphics[scale=0.23]{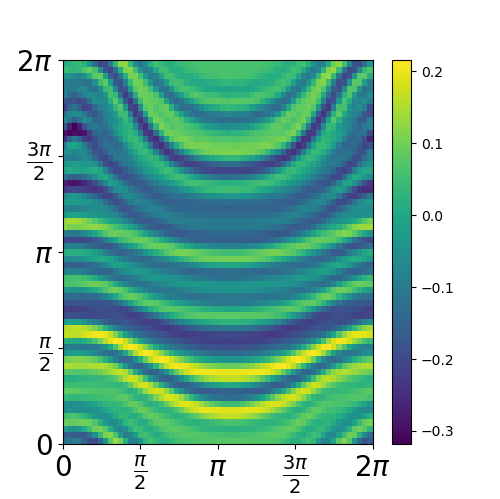}\newline $d=30$
    \end{tabularx}
    \caption{Real part of $\hat{w}_{s,1}q_{y,d}$ for $y=0$ and $s=0.1$.}\label{fig:gaussian}
\end{figure}

We first visualize $\hat{w}_{s,j}$, the eigenvector in Theorem~\ref{prop:hat_w_conti}.
Let $\mcl{Y}=\mathbb{T}$ and $\mcl{Z}=\mathbb{T}^2$.
Consider the dynamical system
\begin{equation}
\bigg(\frac{\mr{d}y(t)}{\mr{d}t},\bigg(\frac{\mr{d}z_1(t)}{\mr{d}t},\frac{\mr{d}z_2(t)}{\mr{d}t}\bigg)\bigg)=\bigg(1,\bigg(-\frac{\partial \zeta}{\partial z_2}(y(t),z(t)),\frac{\partial \zeta}{\partial z_1}(y(t),z(t))\bigg)\bigg),\label{eq:ex_stream_fun}
\end{equation}
where $\zeta(y,z)=\mr{e}^{\kappa(\cos(z_1-y)+\cos z_2)}$.
This problem is also studied by Giannakis and Das~\cite{giannakis20}.
In this case, $\frac{\partial g}{\partial t}(0,y,z)=(-\frac{\partial \zeta}{\partial z_2}(y,z),\frac{\partial \zeta}{\partial z_1}(y,z))$ and $V_{\Phi}=\frac{\partial}{\partial y}-\frac{\partial \zeta}{\partial z_2}(y,z)\frac{\partial}{\partial z_1}+\frac{\partial \zeta}{\partial z_1}(y,z)\frac{\partial}{\partial z_2}$.
To construct the subspace $\mcl{V}_j(y)$ approximately, we set $\kappa=0.5$ and approximated $V_{\Phi}$ in the RKHS $\hil_{\tau}\otimes \hil_{\tau}\otimes \hil_{\tau}$ with $\tau=0.1$.
Here, $\hil_{\tau}$ is the RKHS associated with the positive definite kernel $p_{\tau}:\mathbb{T}\times\mathbb{T}\to\mathbb{C}$ defined as $p_{\tau}(y_1,y_2)=\sum_{i\in\mathbb{Z}}\lambda_{\tau,i}\phi_i(y_1)\overline{\phi_i(y_2)}$.
In addition, $\phi_i(y)=\mr{e}^{\sqrt{-1}iy}$ and $\lambda_{\tau,i}=\mr{e}^{-\tau\vert i\vert^p}$ for $i\in\mathbb{Z}$.
We set $p=0.1$.
We computed eigenvectors $\tilde{v}_1,\tilde{v}_2,\ldots$ of the approximated operator in the RKHS.
Here, the index is ordered from the eigenvector corresponding to the closest eigenvalue to $10^{-10}$.
We set $\mcl{V}_1(y)=\opn{Span}\{\tilde{v}_1(y,\cdot),\ldots,\tilde{v}_d(y,\cdot)\}$.
Since the eigenvector $w_{s,j}(y)$ is an operator, its visualization is not easy.
Thus, we set any test vector $q_{y,d} \in\mcl{V}_1(y)$ and visualize $w_{s,j}(y)q_{y,d}$ instead of $w_{s,j}$ itself.
As the test vector, we set $q_{y,d}=1/d\sum_{i=1}^d\tilde{v}_i(y,\cdot)$ as an example.
Figure~\ref{fig:gaussian} shows the eigenvector $\hat{w}_{s,1}$ acting on the vector $q_{y,d}$.
We can see that the pattern becomes more clear as $d$ becomes large, which implies that considering higher dimensional subspaces $\mcl{V}_1(y)$ than a one dimensional space catches the feature of the dynamical system in this case.
\newcolumntype{Y}{>{\centering\arraybackslash}X}

\subsection{Idealized stratospheric flow}
\begin{figure}[t]
\begin{tabularx}{\textwidth}{YYY}
\includegraphics[scale=0.25]{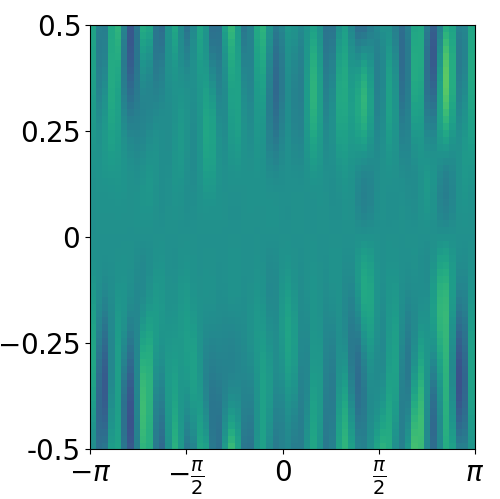}\newline $\sigma(\hat{N}_{s,j})=-11.532\mr{i}$, \newline $\sigma(V_{\Phi})=0.040001\mr{i}$&
\includegraphics[scale=0.25]{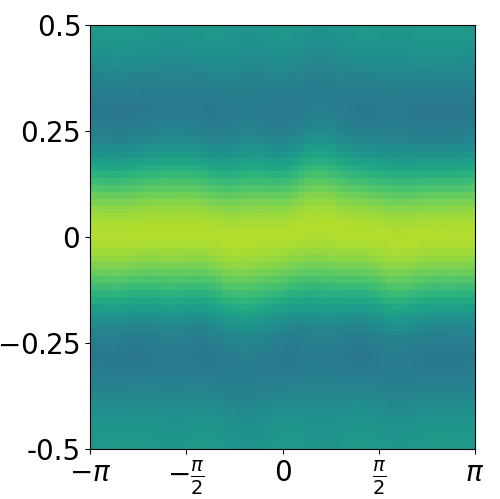}\newline $\sigma(\hat{N}_{s,j})=12.4976\mr{i}$, \newline $\sigma(V_{\Phi})=0.22200\mr{i}$&
\includegraphics[scale=0.25]{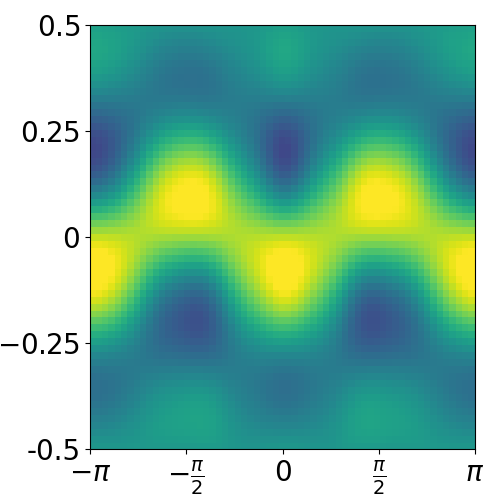}\newline $\sigma(\hat{N}_{s,j})=13.028\mr{i}$, \newline $\sigma(V_{\Phi})=0.17188\mr{i}$
\end{tabularx}
\caption{Real part of $\hat{w}_{s,j}\tilde{v}_j(y,\cdot)$ for $s=y=0$}\label{fig:jet}
\end{figure}

Next, we observe the eigenoperator $\hat{N}_{s,j}$. 
We study what information the eigenoperators capture.
Let $\mcl{Y}=\mathbb{T}$ and $\mcl{Z}=\mathbb{T}\times [-\pi,\pi]$.
Consider the same dynamical system as Eq.~\eqref{eq:ex_stream_fun} with $\zeta(y,z)=c_3z_2-U_0L\tanh(z_2/L)+\sum_{i=1}^3A_iU_0L\opn{sech}^2(z_2/L)\cos(k_iz_1-\sigma_iy)$, where $L=0.1$, $A_1=0.075$, $A_2=0.4$, $A_3=0.2$, $k_1=1$, $k_2=2k_1$, $k_3=3k_1$, $U_0=62.66$, $c_3=0.7U_0$, $\sigma_2=-1$, and $\sigma_1=2\sigma_1$.
A similar problem is also studied by Froyland et al.~\cite{froyland10b}.
We approximated $V_{\Phi}$ in the RKHS $\hil_{\tau}\otimes \hil_{\tau}\otimes \hil_{\tau}$ with $\tau=0.1$.
We computed eigenvectors $\tilde{v}_1,\tilde{v}_2,\ldots$ of the approximated operator in the RKHS and set $\mcl{V}_j(y)=\opn{Span}\{\tilde{v}_j(y)\}$ for $j=1,2,\ldots$.
We study $\hat{N}_{s,j}$ for different $j$ and what information we can extract according to $\hat{N}_{s,j}$.
Figure~\ref{fig:jet} shows the heatmap of the function $\hat{w}_{0,j}\tilde{v}_j(0,\cdot)$ for different values of $j$.
Since we have $\hat{w}_{0,j}\tilde{v}_j(y,\cdot)=\tilde{v}_j(y,\cdot)$, it provides us coherent patterns.
We computed the spectrum of $\hat{N}_{s,j}$ for the corresponding $j$.
The eigenoperator $\hat{N}_{s,j}$ is different from the spectrum of $V_{\Phi}$, the generator of the Koopman operator.
We also computed the spectrum of $V_{\Phi}$.
We can see that the pattern becomes complicated as the magnitude of the spectrum $\sigma(\hat{N}_{0,j})$ of the eigenoperator becomes large.
On the other hand, the spectrum of $V_{\Phi}$ does not provide such an observation.

\section{Conclusion and discussion}
\label{sec:conclusions}

In this paper, we considered a skew product dynamical system on $\mcl{Y}\times \mcl{Z}$ and defined a linear operator on a Hilbert $C^*$-module related to the Koopman operator.
We proposed the eigenoperator decomposition as a generalization of the eigenvalue decomposition.
The eigenvectors are constructed using a cocycle.
The eigenoperators reconstruct the Koopman operator projected on generalized Oseledets subspaces.
Thus, if the Oseledets subspaces are infinite-dimensional spaces, the eigenoperators can have continuous spectra related to the Koopman operator.
Our approach is different from existing approach to dealing with continuous and residual spectra of Koopman operators, such as focusing on the spectral measure~\cite{korda20} and approximating Koopman operators using compact operators on different space from the space where the Koopman operators are defined~\cite{slipantschuk20,giannakis20}.
In addition, the proposed decomposition gives us information of the behavior of coherent patterns on $\mcl{Z}$.
Extracting coherent structure of skew product dynamical systems has been investigated~\cite{FROYLAND10a,froyland10b,das21}.
The proposed decomposition will allow us to classify these coherent patterns.

For future work, studying data-driven approaches to obtaining the decomposition is an important direction of researches.
Investigating practical and computationally efficient ways to approximate operators on Hilbert $C^*$-modules would be essential in that direction of researches.
Another interesting direction is applying the proposed decomposition to quantum computation.
Decomposing a Koopman operator for quantum computation was proposed~\cite{giannakis22}.
It would be interesting to generalize the decomposition using the proposed decomposition.

\section*{Acknowledgments}
We thank Suddhasattwa Das for many constructive discussions on this work. DG acknowledges support from the US National Science Foundation under grant DMS-1854383, the US Office of Naval Research under MURI grant N00014-19-1-242, and the US Department of Defense, Basic Research Office under Vannevar Bush Faculty Fellowship grant N00014-21-1-2946.
MI and II acknowledge support from  the Japan Science and Technlogy Agency under CREST grant JPMJCR1913.
II acknowledge support from  the Japan Science and Technlogy Agency under ACT-X grant JPMJAX2004.

\begin{appendices}
\section{Proof of Proposition~\ref{prop:K_TinM}}\label{ap:proof_propK_T}
To show Proposition~\ref{prop:K_TinM}, we use the following lemmas~\cite{skeide00}.
\begin{lemma}\label{lem:linear_op}
We have $\modu\subseteq \mcl{B}(\lz,\lx)$.
\end{lemma}
\begin{proof}
Let $\iota:\ly\otimes_{\opn{alg}} \alg\to\mcl{B}(\lz,\lx)$ be defined as
$(\iota(v\otimes a)u)(y,z)=v(y)(au)(z)$.
Then, $\iota$ is an injection.
In addition, we have
\begin{align*}
\Vert \iota(v\otimes a)\Vert_{\mcl{B}(\lz,\lx)}
&=\sup_{\Vert u\Vert_{\lz}=1}\Vert \iota(v\otimes a)u\Vert_{\lx}\\
&=\Vert v\Vert_{\ly}\sup_{\Vert u\Vert_{\lz}=1}\Vert au\Vert_{\lz}
=\Vert v\otimes a\Vert_{\modu}.
\end{align*}
Therefore, we have $\modu\subseteq \mcl{B}(\lz,\lx)$.
\end{proof}
We recall that a Hilbert $C^*$-module $\modu$ is referred to as self-dual if for any bounded $\alg$-linear map $b:\modu\to\alg$, there exists a unique $\hat{b}\in\modu$ such that $b(w)=\sbracket{\hat{b},w}_{\modu}$.
\begin{lemma}\label{lem:strongly_closed}
Let $\{w_i\}_{i=1}^{\infty}$ be a sequence in $\modu$.
Assume there exists $w\in\mcl{B}(\lz,\lx)$ such that for any $u\in\lz$, $\iota(w_i)u\to wu$ in $\lx$, where $\iota$ is defined in the proof of Lemma~\ref{lem:linear_op}.
Then, $w\in\modu$.
\end{lemma}
\begin{proof}
For $\tilde{w}\in\modu$ and $u_1,u_2\in\lz$, we have
\begin{align*}
\bracket{\iota(\tilde{w})u_1,wu_2}=\bracket{\iota(\tilde{w})u_1,\lim_{i\to\infty}\iota(w_i)u_2}
&=\lim_{i\to\infty}\int_{z\in\mcl{Z}}\int_{y\in\mcl{Y}}\overline{u_1(z)}\tilde{w}(y)^*w_i(y)u_2(z)\mr{d}\mu(y)\mr{d}\nu(z)\\
&=\lim_{i\to\infty}\bracket{u_1,\bracket{\tilde{w},w_i}_{\modu}u_2}_{\lz}.
\end{align*}
Thus, by the Riesz representation theorem, there exists $a\in\alg$ such that $\bracket{\iota(\tilde{w})u_1,wu_2}=\bracket{u_1,au_2}$. 
The map $\tilde{w}\mapsto a^*$ is a bounded $\alg$-linear map from $\modu$ to $\alg$.
Since $\alg$ is self-dual, $\modu$ is also self-dual. 
As a result, there exists $\hat{w}\in\modu$ such that $\sbracket{\hat{w},\tilde{w}}_{\modu}=a^*$ and $w=\iota(\hat{w})$.
\end{proof}
\begin{lemma}\label{lem:winM}
Let $w:\mcl{Y}\to\alg$.
Assume for any $u\in\lz$, the map $(y,z)\mapsto (w(y)u)(z)$ is contained in $\lx$.
Then, $w\in\modu$.
\end{lemma}
\begin{proof}
Let $\{\gamma_i\}_{i=1}^{\infty}$ be an orthonormal basis of $\lz$.
For $u\in\ly$, we have
\begin{equation*}
w(y)u=w(y)\sum_{i=1}^{\infty}\gamma_i\gamma_i'u=\sum_{i=1}^{\infty}w(y)\gamma_i\gamma_i'u.
\end{equation*}
Here, $\gamma_i'$ denotes the dual of $\gamma_i\in \lz$.
Since the map $(y,z)\mapsto (w(y)\gamma_i)(z)$ is in $\lx$, we obtain $w(y)\gamma_i\gamma_i'\in\modu$.
Thus, by regarding $w$ as an element in $\mcl{B}(\lz,\lx)$ defined as $u\mapsto ((y,z)\mapsto (w(y)u)(z))$, by Lemma~\ref{lem:strongly_closed}, $w\in\modu$ holds.
\end{proof}
\begin{proof}[Proof of Proposition~\ref{prop:K_TinM}]
Since $(y,z)\mapsto (K_T(v\otimes a)(y)u)(z)=v(h(y))(au)(g(y,z))$ is in $\lx$,
by Lemma~\ref{lem:winM}, $K_T(v\otimes a)\in\modu$ holds.

Regarding the unitarity of $K_T$, let $L_T:\modu\to\modu$ be a right $\alg$-linear operator defined as $L_T(v\otimes a)=v(h^{-1}(y))U_{g(h^{-1}(y),\cdot)}^*a$.
Then, $L_T$ is the inverse of $K_T$ and for $v_1,v_2\in\ly$ and $a_1,a_2\in\alg$, we have
\begin{align*}
&\bracket{K_T(v_1\otimes a_1),K_T(v_2\otimes a_2)}_{\modu}= \int_{y\in\mcl{Y}}\overline{v_1(h(y))}a_1^*U_{g(y,\cdot)}^*U_{g(y,\cdot)}a_2v_2(h(y))\mr{d}\mu(y)\\
&\qquad\qquad=\bracket{v_1,v_2}_{\ly}a_1^*a_2=\bracket{v_1\otimes a_1,v_2\otimes a_2}_{\modu}.
\end{align*}
\end{proof}

\end{appendices}



\begin{thebibliography}{44}
\ifx \bisbn   \undefined \def \bisbn  #1{ISBN #1}\fi
\ifx \binits  \undefined \def \binits#1{#1}\fi
\ifx \bauthor  \undefined \def \bauthor#1{#1}\fi
\ifx \batitle  \undefined \def \batitle#1{#1}\fi
\ifx \bjtitle  \undefined \def \bjtitle#1{#1}\fi
\ifx \bvolume  \undefined \def \bvolume#1{\textbf{#1}}\fi
\ifx \byear  \undefined \def \byear#1{#1}\fi
\ifx \bissue  \undefined \def \bissue#1{#1}\fi
\ifx \bfpage  \undefined \def \bfpage#1{#1}\fi
\ifx \blpage  \undefined \def \blpage #1{#1}\fi
\ifx \burl  \undefined \def \burl#1{\textsf{#1}}\fi
\ifx \doiurl  \undefined \def \doiurl#1{\url{https://doi.org/#1}}\fi
\ifx \betal  \undefined \def \betal{\textit{et al.}}\fi
\ifx \binstitute  \undefined \def \binstitute#1{#1}\fi
\ifx \binstitutionaled  \undefined \def \binstitutionaled#1{#1}\fi
\ifx \bctitle  \undefined \def \bctitle#1{#1}\fi
\ifx \beditor  \undefined \def \beditor#1{#1}\fi
\ifx \bpublisher  \undefined \def \bpublisher#1{#1}\fi
\ifx \bbtitle  \undefined \def \bbtitle#1{#1}\fi
\ifx \bedition  \undefined \def \bedition#1{#1}\fi
\ifx \bseriesno  \undefined \def \bseriesno#1{#1}\fi
\ifx \blocation  \undefined \def \blocation#1{#1}\fi
\ifx \bsertitle  \undefined \def \bsertitle#1{#1}\fi
\ifx \bsnm \undefined \def \bsnm#1{#1}\fi
\ifx \bsuffix \undefined \def \bsuffix#1{#1}\fi
\ifx \bparticle \undefined \def \bparticle#1{#1}\fi
\ifx \barticle \undefined \def \barticle#1{#1}\fi
\bibcommenthead
\ifx \bconfdate \undefined \def \bconfdate #1{#1}\fi
\ifx \botherref \undefined \def \botherref #1{#1}\fi
\ifx \url \undefined \def \url#1{\textsf{#1}}\fi
\ifx \bchapter \undefined \def \bchapter#1{#1}\fi
\ifx \bbook \undefined \def \bbook#1{#1}\fi
\ifx \bcomment \undefined \def \bcomment#1{#1}\fi
\ifx \oauthor \undefined \def \oauthor#1{#1}\fi
\ifx \citeauthoryear \undefined \def \citeauthoryear#1{#1}\fi
\ifx \endbibitem  \undefined \def \endbibitem {}\fi
\ifx \bconflocation  \undefined \def \bconflocation#1{#1}\fi
\ifx \arxivurl  \undefined \def \arxivurl#1{\textsf{#1}}\fi
\csname PreBibitemsHook\endcsname

\bibitem[\protect\citeauthoryear{Koopman}{1931}]{koopman31}
\begin{barticle}
\bauthor{\bsnm{Koopman}, \binits{B.}}:
\batitle{Hamiltonian systems and transformation in {H}ilbert space}.
\bjtitle{Proc. Natl. Acad. Sci.}
\bvolume{17}(\bissue{5}),
\bfpage{315}--\blpage{318}
(\byear{1931})
\doiurl{10.1073/pnas.17.5.315}
\end{barticle}
\endbibitem

\bibitem[\protect\citeauthoryear{Koopman and von
  Neumann}{1932}]{KoopmanVonNeumann32}
\begin{barticle}
\bauthor{\bsnm{Koopman}, \binits{B.O.}},
\bauthor{\bsnm{Neumann}, \binits{J.}}:
\batitle{Dynamical systems of continuous spectra}.
\bjtitle{Proc. Natl. Acad. Sci.}
\bvolume{18}(\bissue{3}),
\bfpage{255}--\blpage{263}
(\byear{1932})
\doiurl{10.1073/pnas.18.3.255}
\end{barticle}
\endbibitem

\bibitem[\protect\citeauthoryear{Baladi}{2000}]{Baladi00}
\begin{bbook}
\bauthor{\bsnm{Baladi}, \binits{V.}}:
\bbtitle{Positive Transfer Operators and Decay of Correlations}.
\bsertitle{Advanced Series in Nonlinear Dynamics},
vol. \bseriesno{16}.
\bpublisher{World Scientific},
\blocation{Singapore}
(\byear{2000})
\end{bbook}
\endbibitem

\bibitem[\protect\citeauthoryear{Eisner et~al.}{2015}]{EisnerEtAl15}
\begin{bbook}
\bauthor{\bsnm{Eisner}, \binits{T.}},
\bauthor{\bsnm{Farkas}, \binits{B.}},
\bauthor{\bsnm{Haase}, \binits{M.}},
\bauthor{\bsnm{Nagel}, \binits{R.}}:
\bbtitle{Operator Theoretic Aspects of Ergodic Theory}.
\bsertitle{Graduate Texts in Mathematics},
vol. \bseriesno{272}.
\bpublisher{Springer},
\blocation{Cham}
(\byear{2015})
\end{bbook}
\endbibitem

\bibitem[\protect\citeauthoryear{Froyland}{1997}]{Froyland97}
\begin{barticle}
\bauthor{\bsnm{Froyland}, \binits{G.}}:
\batitle{Computer-assisted bounds for the rate of decay of correlations}.
\bjtitle{Commun. Math. Phys.}
\bvolume{189}(\bissue{NN}),
\bfpage{237}--\blpage{257}
(\byear{1997})
\doiurl{10.1007/s002200050198}
\end{barticle}
\endbibitem

\bibitem[\protect\citeauthoryear{Dellnitz and Junge}{1999}]{DellnitzJunge99}
\begin{barticle}
\bauthor{\bsnm{Dellnitz}, \binits{M.}},
\bauthor{\bsnm{Junge}, \binits{O.}}:
\batitle{On the approximation of complicated dynamical behavior}.
\bjtitle{SIAM J. Numer. Anal.}
\bvolume{36},
\bfpage{491}
(\byear{1999})
\doiurl{10.1137/S0036142996313002}
\end{barticle}
\endbibitem

\bibitem[\protect\citeauthoryear{Dellnitz et~al.}{2000}]{DellnitzEtAl00}
\begin{barticle}
\bauthor{\bsnm{Dellnitz}, \binits{M.}},
\bauthor{\bsnm{Froyland}, \binits{G.}},
\bauthor{\bsnm{Sertl}, \binits{S.}}:
\batitle{On the isolated spectrum of the {P}erron–{F}robenius operator}.
\bjtitle{Nonlinearity}
\bvolume{13},
\bfpage{1171}--\blpage{1188}
(\byear{2000})
\doiurl{10.1088/0951-7715/13/4/310}
\end{barticle}
\endbibitem

\bibitem[\protect\citeauthoryear{Mezi\'c}{2005}]{Mezic05}
\begin{barticle}
\bauthor{\bsnm{Mezi\'c}, \binits{I.}}:
\batitle{Spectral properties of dynamical systems, model reduction and
  decompositions}.
\bjtitle{Nonlinear Dyn.}
\bvolume{41},
\bfpage{309}--\blpage{325}
(\byear{2005})
\doiurl{10.1007/s11071-005-2824-x}
\end{barticle}
\endbibitem

\bibitem[\protect\citeauthoryear{Giannakis et~al.}{2015}]{GiannakisEtAl15}
\begin{bchapter}
\bauthor{\bsnm{Giannakis}, \binits{D.}},
\bauthor{\bsnm{Slawinska}, \binits{J.}},
\bauthor{\bsnm{Zhao}, \binits{Z.}}:
\bctitle{Spatiotemporal feature extraction with data-driven {K}oopman
  operators}.
In: \beditor{\bsnm{Storcheus}, \binits{D.}},
\beditor{\bsnm{Rostamizadeh}, \binits{A.}},
\beditor{\bsnm{Kumar}, \binits{S.}} (eds.)
\bbtitle{Proceedings of the 1st International Workshop on Feature Extraction:
  Modern Questions and Challenges at NIPS 2015}.
\bsertitle{Proceedings of Machine Learning Research},
vol. \bseriesno{44},
pp. \bfpage{103}--\blpage{115}.
\bpublisher{PMLR},
\blocation{Montreal, Canada}
(\byear{2015}).
\burl{https://proceedings.mlr.press/v44/giannakis15.html}
\end{bchapter}
\endbibitem

\bibitem[\protect\citeauthoryear{Klus et~al.}{2020}]{klus17}
\begin{barticle}
\bauthor{\bsnm{Klus}, \binits{S.}},
\bauthor{\bsnm{Schuster}, \binits{I.}},
\bauthor{\bsnm{Muandet}, \binits{K.}}:
\batitle{Eigendecompositions of transfer operators in reproducing kernel
  {H}ilbert spaces}.
\bjtitle{J. Nonlinear Sci.}
\bvolume{30},
\bfpage{283}--\blpage{315}
(\byear{2020})
\doiurl{10.1007/s00332-019-09574-z}
\end{barticle}
\endbibitem

\bibitem[\protect\citeauthoryear{Ishikawa et~al.}{2018}]{ishikawa18}
\begin{bchapter}
\bauthor{\bsnm{Ishikawa}, \binits{I.}},
\bauthor{\bsnm{Fujii}, \binits{K.}},
\bauthor{\bsnm{Ikeda}, \binits{M.}},
\bauthor{\bsnm{Hashimoto}, \binits{Y.}},
\bauthor{\bsnm{Kawahara}, \binits{Y.}}:
\bctitle{Metric on nonlinear dynamical systems with {P}erron-{F}robenius
  operators}.
In: \bbtitle{Proceedings of Advances in Neural Information Processing Systems
  32 (NeurIPS)}
(\byear{2018})
\end{bchapter}
\endbibitem

\bibitem[\protect\citeauthoryear{Hashimoto et~al.}{2020}]{hashimoto20}
\begin{barticle}
\bauthor{\bsnm{Hashimoto}, \binits{Y.}},
\bauthor{\bsnm{Ishikawa}, \binits{I.}},
\bauthor{\bsnm{Ikeda}, \binits{M.}},
\bauthor{\bsnm{Matsuo}, \binits{Y.}},
\bauthor{\bsnm{Kawahara}, \binits{Y.}}:
\batitle{Krylov subspace method for nonlinear dynamical systems with random
  noise}.
\bjtitle{J. Mach. Learn. Res.}
\bvolume{21}(\bissue{172}),
\bfpage{1}--\blpage{29}
(\byear{2020})
\end{barticle}
\endbibitem

\bibitem[\protect\citeauthoryear{Schmid}{2010}]{Schmid10}
\begin{barticle}
\bauthor{\bsnm{Schmid}, \binits{P.J.}}:
\batitle{Dynamic mode decomposition of numerical and experimental data}.
\bjtitle{J. Fluid Mech.}
\bvolume{656},
\bfpage{5}--\blpage{28}
(\byear{2010})
\doiurl{10.1017/S0022112010001217}
\end{barticle}
\endbibitem

\bibitem[\protect\citeauthoryear{Rowley et~al.}{2009}]{RowleyEtAl09}
\begin{barticle}
\bauthor{\bsnm{Rowley}, \binits{C.W.}},
\bauthor{\bsnm{Mezi\'c}, \binits{I.}},
\bauthor{\bsnm{Bagheri}, \binits{S.}},
\bauthor{\bsnm{Schlatter}, \binits{P.}},
\bauthor{\bsnm{Henningson}, \binits{D.S.}}:
\batitle{Spectral analysis of nonlinear flows}.
\bjtitle{J. Fluid Mech.}
\bvolume{641},
\bfpage{115}--\blpage{127}
(\byear{2009})
\doiurl{10.1017/s0022112009992059}
\end{barticle}
\endbibitem

\bibitem[\protect\citeauthoryear{Williams et~al.}{2015}]{williams15}
\begin{barticle}
\bauthor{\bsnm{Williams}, \binits{M.O.}},
\bauthor{\bsnm{Kevrekidis}, \binits{I.G.}},
\bauthor{\bsnm{Rowley}, \binits{C.W.}}:
\batitle{A data–driven approximation of the {K}oopman operator: extending
  dynamic mode decomposition}.
\bjtitle{J. Nonlinear Sci.}
\bvolume{25},
\bfpage{1307}--\blpage{1346}
(\byear{2015})
\doiurl{10.1007/s00332-015-9258-5}
\end{barticle}
\endbibitem

\bibitem[\protect\citeauthoryear{Kawahara}{2016}]{kawahara16}
\begin{bchapter}
\bauthor{\bsnm{Kawahara}, \binits{Y.}}:
\bctitle{Dynamic mode decomposition with reproducing kernels for {K}oopman
  spectral analysis}.
In: \bbtitle{Proceedings of Advances in Neural Information Processing Systems
  30 (NIPS)}
(\byear{2016})
\end{bchapter}
\endbibitem

\bibitem[\protect\citeauthoryear{Arbabi and Mezi\'{c}}{2017}]{hassan17}
\begin{barticle}
\bauthor{\bsnm{Arbabi}, \binits{H.}},
\bauthor{\bsnm{Mezi\'{c}}, \binits{I.}}:
\batitle{Ergodic theory, dynamic mode decomposition, and computation of
  spectral properties of the {K}oopman operator}.
\bjtitle{SIAM J. Appl. Dyn. Syst.}
\bvolume{16}(\bissue{4}),
\bfpage{2096}--\blpage{2126}
(\byear{2017})
\doiurl{10.1137/17M1125}
\end{barticle}
\endbibitem

\bibitem[\protect\citeauthoryear{Rosenfeld et~al.}{2022}]{RosenfeldEtAl22}
\begin{botherref}
\oauthor{\bsnm{Rosenfeld}, \binits{J.A.}},
\oauthor{\bsnm{Kamalapurkar}, \binits{R.}},
\oauthor{\bsnm{Gruss}, \binits{L.F.}},
\oauthor{\bsnm{Johnson}, \binits{T.T.}}:
Dynamic mode decomposition for continuous time systems with the {L}iouville
  operator.
J. Nonlinear Sci.
\textbf{32}
(2022)
\doiurl{10.1007/s00332-021-09746-w}
\end{botherref}
\endbibitem

\bibitem[\protect\citeauthoryear{Korda et~al.}{2020}]{korda20}
\begin{barticle}
\bauthor{\bsnm{Korda}, \binits{M.}},
\bauthor{\bsnm{Putinar}, \binits{M.}},
\bauthor{\bsnm{Mezi\'{c}}, \binits{I.}}:
\batitle{Data-driven spectral analysis of the {K}oopman operator}.
\bjtitle{Appl. Comput. Harmon. Anal.}
\bvolume{48}(\bissue{2}),
\bfpage{599}--\blpage{629}
(\byear{2020})
\doiurl{10.1016/j.acha.2018.08.002}
\end{barticle}
\endbibitem

\bibitem[\protect\citeauthoryear{Slipantschuk et~al.}{2020}]{slipantschuk20}
\begin{barticle}
\bauthor{\bsnm{Slipantschuk}, \binits{J.}},
\bauthor{\bsnm{Bandtlow}, \binits{O.F.}},
\bauthor{\bsnm{Just}, \binits{W.}}:
\batitle{Dynamic mode decomposition for analytic maps}.
\bjtitle{Commun. Nonlinear Sci. Numer. Simul.}
\bvolume{84},
\bfpage{105179}
(\byear{2020})
\doiurl{10.1016/j.cnsns.2020.105179}
\end{barticle}
\endbibitem

\bibitem[\protect\citeauthoryear{Colbrook and
  Townsend}{2021}]{ColbrookTownsend21}
\begin{botherref}
\oauthor{\bsnm{Colbrook}, \binits{M.J.}},
\oauthor{\bsnm{Townsend}, \binits{A.}}:
Rigorous Data-Driven Computation of Spectral Properties of {K}oopman Operators
  for Dynamical Systems
(2021).
\url{https://arxiv.org/abs/2111.14889}
\end{botherref}
\endbibitem

\bibitem[\protect\citeauthoryear{Das et~al.}{2021}]{das21}
\begin{barticle}
\bauthor{\bsnm{Das}, \binits{S.}},
\bauthor{\bsnm{Giannakis}, \binits{D.}},
\bauthor{\bsnm{Slawinska}, \binits{J.}}:
\batitle{Reproducing kernel {H}ilbert space compactification of unitary
  evolution groups}.
\bjtitle{Appl. Comput. Harmon. Anal.}
\bvolume{54},
\bfpage{75}--\blpage{136}
(\byear{2021})
\doiurl{10.1016/j.acha.2021.02.004}
\end{barticle}
\endbibitem

\bibitem[\protect\citeauthoryear{Ulam}{1964}]{Ulam64}
\begin{bbook}
\bauthor{\bsnm{Ulam}, \binits{S.M.}}:
\bbtitle{Problems in Modern Mathematics}.
\bpublisher{Dover Publications},
\blocation{Mineola}
(\byear{1964})
\end{bbook}
\endbibitem

\bibitem[\protect\citeauthoryear{Blank et~al.}{2002}]{BlankEtAl02}
\begin{barticle}
\bauthor{\bsnm{Blank}, \binits{M.}},
\bauthor{\bsnm{Keller}, \binits{G.}},
\bauthor{\bsnm{Liverani}, \binits{C.}}:
\batitle{Ruelle--{P}erron--{F}robenius spectrum for {A}nosov maps}.
\bjtitle{Nonlinearity}
\bvolume{15}(\bissue{6}),
\bfpage{1905}--\blpage{1973}
(\byear{2002})
\doiurl{10.1088/0951-7715/15/6/309}
\end{barticle}
\endbibitem

\bibitem[\protect\citeauthoryear{Junge and Koltai}{2009}]{JungeKoltai09}
\begin{barticle}
\bauthor{\bsnm{Junge}, \binits{O.}},
\bauthor{\bsnm{Koltai}, \binits{P.}}:
\batitle{Discretization of the {F}robenius--{P}erron operator using a sparse
  {H}aar tensor basis: {T}he sparse {U}lam method}.
\bjtitle{SIAM J. Numer. Anal.}
\bvolume{47},
\bfpage{3464}--\blpage{2485}
(\byear{2009})
\doiurl{10.1137/080716864}
\end{barticle}
\endbibitem

\bibitem[\protect\citeauthoryear{Froyland et~al.}{2010}]{froyland10b}
\begin{barticle}
\bauthor{\bsnm{Froyland}, \binits{G.}},
\bauthor{\bsnm{Santitissadeekorn}, \binits{N.}},
\bauthor{\bsnm{Monahan}, \binits{A.}}:
\batitle{Transport in time-dependent dynamical systems: Finite-time coherent
  sets}.
\bjtitle{Chaos: An Interdisciplinary Journal of Nonlinear Science}
\bvolume{20}(\bissue{4}),
\bfpage{043116}
(\byear{2010})
\doiurl{10.1063/1.3502450}
\end{barticle}
\endbibitem

\bibitem[\protect\citeauthoryear{Froyland and Koltai}{2017}]{froyland17}
\begin{barticle}
\bauthor{\bsnm{Froyland}, \binits{G.}},
\bauthor{\bsnm{Koltai}, \binits{P.}}:
\batitle{Estimating long-term behavior of periodically driven flows without
  trajectory integration}.
\bjtitle{Nonlinearity}
\bvolume{30}(\bissue{5}),
\bfpage{1948}
(\byear{2017})
\doiurl{10.1088/1361-6544/aa6693}
\end{barticle}
\endbibitem

\bibitem[\protect\citeauthoryear{Giannakis and Das}{2020}]{giannakis20}
\begin{barticle}
\bauthor{\bsnm{Giannakis}, \binits{D.}},
\bauthor{\bsnm{Das}, \binits{S.}}:
\batitle{Extraction and prediction of coherent patterns in incompressible flows
  through space-time {K}oopman analysis}.
\bjtitle{Physica D: Nonlinear Phenomena}
\bvolume{402},
\bfpage{132211}
(\byear{2020})
\doiurl{10.1016/j.physd.2019.132211}
\end{barticle}
\endbibitem

\bibitem[\protect\citeauthoryear{Oseledets}{1968}]{Oseledets68}
\begin{barticle}
\bauthor{\bsnm{Oseledets}, \binits{V.I.}}:
\batitle{A multiplicative ergodic theorem}.
\bjtitle{Trans. Moscow Math. Soc.}
\bvolume{19},
\bfpage{197}--\blpage{231}
(\byear{1968})
\end{barticle}
\endbibitem

\bibitem[\protect\citeauthoryear{Ruelle}{1968}]{Ruelle68}
\begin{barticle}
\bauthor{\bsnm{Ruelle}, \binits{D.}}:
\batitle{Statistical mechanics of a one-dimensional lattice gas}.
\bjtitle{Comm. Math. Phys.}
\bvolume{9}(\bissue{4}),
\bfpage{267}--\blpage{278}
(\byear{1968})
\end{barticle}
\endbibitem

\bibitem[\protect\citeauthoryear{Thieullen}{1987}]{thieullen87}
\begin{barticle}
\bauthor{\bsnm{Thieullen}, \binits{P.}}:
\batitle{Fibres dynamiques asymptotiquement compacts exposants de lyapounov.
  entropie. dimension}.
\bjtitle{Annales de l'Institut Henri Poincar\'{e} C, Analyse non lin\'{e}aire}
\bvolume{4}(\bissue{1}),
\bfpage{49}--\blpage{97}
(\byear{1987})
\doiurl{10.1016/S0294-1449(16)30373-0}
\end{barticle}
\endbibitem

\bibitem[\protect\citeauthoryear{Schauml\"{o}ffel}{1991}]{Schaumloffel91}
\begin{bchapter}
\bauthor{\bsnm{Schauml\"{o}ffel}, \binits{K.}}:
\bctitle{Multiplicative ergodic theorems in infinite dimensions}.
In: \bbtitle{Lyapunov Exponents. Lecture Notes in Mathematics},
vol. \bseriesno{1486}.
\bpublisher{Springer},
\blocation{Heidelberg}
(\byear{1991})
\end{bchapter}
\endbibitem

\bibitem[\protect\citeauthoryear{Froyland et~al.}{2010}]{FroylandEtAl10c}
\begin{barticle}
\bauthor{\bsnm{Froyland}, \binits{G.}},
\bauthor{\bsnm{Lloyd}, \binits{S.}},
\bauthor{\bsnm{Quas}, \binits{A.}}:
\batitle{Coherent structures and isolated spectrum for {P}erron--{F}robenius
  cocycles}.
\bjtitle{Ergod. Theory Dyn. Syst.}
\bvolume{30}(\bissue{3}),
\bfpage{729}--\blpage{756}
(\byear{2010})
\doiurl{10.1017/S0143385709000339}
\end{barticle}
\endbibitem

\bibitem[\protect\citeauthoryear{Gonz\'{a}lez-Tokman and
  Quas}{2014}]{gonzalez14}
\begin{barticle}
\bauthor{\bsnm{Gonz\'{a}lez-Tokman}, \binits{C.}},
\bauthor{\bsnm{Quas}, \binits{A.}}:
\batitle{A semi-invertible operator oseledets theorem}.
\bjtitle{Ergod. Theory Dyn. Syst.}
\bvolume{34}(\bissue{4}),
\bfpage{1230}--\blpage{1272}
(\byear{2014})
\doiurl{10.1017/etds.2012.189}
\end{barticle}
\endbibitem

\bibitem[\protect\citeauthoryear{Froyland et~al.}{2010}]{FROYLAND10a}
\begin{barticle}
\bauthor{\bsnm{Froyland}, \binits{G.}},
\bauthor{\bsnm{Lloyd}, \binits{S.}},
\bauthor{\bsnm{Santitissadeekorn}, \binits{N.}}:
\batitle{Coherent sets for nonautonomous dynamical systems}.
\bjtitle{Physica D: Nonlinear Phenomena}
\bvolume{239}(\bissue{16}),
\bfpage{1527}--\blpage{1541}
(\byear{2010})
\doiurl{10.1016/j.physd.2010.03.009}
\end{barticle}
\endbibitem

\bibitem[\protect\citeauthoryear{Froyland and Junge}{2018}]{froyland18}
\begin{barticle}
\bauthor{\bsnm{Froyland}, \binits{G.}},
\bauthor{\bsnm{Junge}, \binits{O.}}:
\batitle{Robust {FEM}-based extraction of finite-time coherent sets using
  scattered, sparse, and incomplete trajectories}.
\bjtitle{SIAM J. Appl. Dyn. Syst.}
\bvolume{17}(\bissue{2}),
\bfpage{1891}--\blpage{1924}
(\byear{2018})
\doiurl{10.1137/17M1129738}
\end{barticle}
\endbibitem

\bibitem[\protect\citeauthoryear{Froyland}{2015}]{froyland15}
\begin{barticle}
\bauthor{\bsnm{Froyland}, \binits{G.}}:
\batitle{Dynamic isoperimetry and the geometry of {L}agrangian coherent
  structures}.
\bjtitle{Nonlinearity}
\bvolume{28}(\bissue{10}),
\bfpage{3587}
(\byear{2015})
\doiurl{10.1088/0951-7715/28/10/3587}
\end{barticle}
\endbibitem

\bibitem[\protect\citeauthoryear{Lance}{1995}]{lance95}
\begin{bbook}
\bauthor{\bsnm{Lance}, \binits{E.C.}}:
\bbtitle{Hilbert {$C^*$}-modules -- a Toolkit for Operator Algebraists}.
\bsertitle{London Mathematical Society Lecture Note Series, vol. 210}.
\bpublisher{Cambridge University Press},
\blocation{New York}
(\byear{1995})
\end{bbook}
\endbibitem

\bibitem[\protect\citeauthoryear{Froyland and Koltai}{2023}]{FroylandKoltai23}
\begin{barticle}
\bauthor{\bsnm{Froyland}, \binits{G.}},
\bauthor{\bsnm{Koltai}, \binits{P.}}:
\batitle{Detecting the birth and death of finite-time coherent sets}.
\bjtitle{Commun. Pure Appl. Math}
(\byear{2023})
\doiurl{10.1002/cpa.22115}
\end{barticle}
\endbibitem

\bibitem[\protect\citeauthoryear{Hashimoto et~al.}{2023}]{hashimoto23}
\begin{bchapter}
\bauthor{\bsnm{Hashimoto}, \binits{Y.}},
\bauthor{\bsnm{Komura}, \binits{F.}},
\bauthor{\bsnm{Ikeda}, \binits{M.}}:
\bctitle{{H}ilbert {$C^*$}-module for analyzing structured data}.
In: \bbtitle{Matrix and Operator Equations}.
\bpublisher{Springer},
\blocation{Cham}
(\byear{2023}).
\doiurl{10.1007/16618_2023_58}
\end{bchapter}
\endbibitem

\bibitem[\protect\citeauthoryear{Choe}{1985}]{choe85}
\begin{barticle}
\bauthor{\bsnm{Choe}, \binits{Y.H.}}:
\batitle{{$C_0$}-semigroups on a locally convex space}.
\bjtitle{J. Math. Anal. Appl.}
\bvolume{106}(\bissue{2}),
\bfpage{293}--\blpage{320}
(\byear{1985})
\doiurl{10.1016/0022-247X(85)90115-5}
\end{barticle}
\endbibitem

\bibitem[\protect\citeauthoryear{Eisner et~al.}{2016}]{eisner16}
\begin{bbook}
\bauthor{\bsnm{Eisner}, \binits{T.}},
\bauthor{\bsnm{Farkas}, \binits{B.}},
\bauthor{\bsnm{Haase}, \binits{M.}},
\bauthor{\bsnm{Nagel}, \binits{R.}}:
\bbtitle{Operator Theoretic Aspects of Ergodic Theory}.
\bpublisher{Springer},
\blocation{Cham}
(\byear{2016})
\end{bbook}
\endbibitem

\bibitem[\protect\citeauthoryear{Giannakis et~al.}{2022}]{giannakis22}
\begin{barticle}
\bauthor{\bsnm{Giannakis}, \binits{D.}},
\bauthor{\bsnm{Ourmazd}, \binits{A.}},
\bauthor{\bsnm{Pfeffer}, \binits{P.}},
\bauthor{\bsnm{Schumacher}, \binits{J.}},
\bauthor{\bsnm{Slawinska}, \binits{J.}}:
\batitle{Embedding classical dynamics in a quantum computer}.
\bjtitle{Phys. Rev. A}
\bvolume{105},
\bfpage{052404}
(\byear{2022})
\doiurl{10.1103/PhysRevA.105.052404}
\end{barticle}
\endbibitem

\bibitem[\protect\citeauthoryear{Skeide}{2000}]{skeide00}
\begin{barticle}
\bauthor{\bsnm{Skeide}, \binits{M.}}:
\batitle{Generalised matrix {$C^{\ast}$}-algebras and representations of
  {H}ilbert modules}.
\bjtitle{Math. Proc. Roy. Irish Acad.}
\bvolume{100A}(\bissue{1}),
\bfpage{11}--\blpage{38}
(\byear{2000})
\end{barticle}
\endbibitem

\end{thebibliography}

\end{document}